\documentclass[a4paper,11pt,reqno]{amsart}

\usepackage{listings}

\usepackage{graphicx}
\usepackage{amsmath}
\usepackage{amssymb}
\usepackage{amsfonts}
\usepackage{epsfig}
\usepackage{epigraph}
\usepackage[hang]{caption}
\usepackage{hhline}
\usepackage{comment}
\usepackage{mathrsfs}

\usepackage{graphics}
\usepackage{float}

\usepackage{mathtools}

\usepackage[T1]{fontenc}

\usepackage[utf8]{inputenc}

\DeclareFixedFont{\ttb}{T1}{txtt}{bx}{n}{8} 
\DeclareFixedFont{\ttm}{T1}{txtt}{m}{n}{8}  

\graphicspath{{eps/}}

 \DeclarePairedDelimiter\abs{\lvert}{\rvert}%

\usepackage{color}
\definecolor{deepblue}{rgb}{0,0,0.85}
\definecolor{deepred}{rgb}{0.6,0,0}
\definecolor{deepgreen}{rgb}{0,0.5,0}
\definecolor{deepbrown}{rgb}{0.88,0.52,0.26}
\definecolor{lavender}{rgb}{0.81,0.41,0.83}

\usepackage{listings}

\newcommand\pythonstyle{\lstset{
language=Python,
basicstyle=\ttm,
morekeywords={str},                
keywordstyle=\ttb\color{deepblue},
emph={MyClass,__init__},          
emphstyle=\ttb\color{deepred},    
stringstyle=\color{deepgreen},
commentstyle = \color{deepgreen},
frame=tb,                         
showstringspaces=false,
classoffset=2,
morekeywords={self},
keywordstyle=\color{deepbrown},
}}

\lstnewenvironment{python}[1][]
{
\pythonstyle
\lstset{#1}
}
{}
\newcommand\pythoninline[1]{{\pythonstyle\lstinline!#1!}}

\makeatletter

\newtheorem{theorem}{Theorem}[section]
\newtheorem*{theorem-non}{Theorem}
\newtheorem*{corollary-non}{Corollary}

\@addtoreset{equation}{section}

\def\Ht{{\rm ht}}

\newtheorem{corollary}[theorem]{Corollary}
\newtheorem{remark}[theorem]{Remark}
\newtheorem{proposition}[theorem]{Proposition}
\newtheorem{lemma}[theorem]{Lemma}

\def\PerfProof{{\it Proof.\ }}

\begin{document}

\title[Extraspecial pairs and calculating structure constants]
 {Extraspecial pairs in the multiply-laced root systems and calculating structure constants}
         \author{Rafael Stekolshchik}

\date{}

\begin{abstract}
  The notions of {\it special} and {\it extraspecial pairs} of roots were introduced by Carter
  for calculating structure constants, \cite{Ca72}.
  Let $\{r, s\}$ be a special pair of roots for which the structure constant $N(r,s)$ is sought,
  and let $\{r_1, s_1\}$ be the extraspecial pair of roots corresponding to $\{r, s\}$.
  Consider the ordered set $\{r_1, r, s, s_1\}$, we will call such a set  a {\it quartet}.
  By studying the different quartets,
  we gain additional insight into the internal structure of the root system.
  It is  shown that for the case $B_n$ we can avoid finding $6$ squares of lengths in the formula 
  for calculating the structure constants. 
  The calculation formula for $B_n$ coincides with the formula for the simply-laced case. For the case $C_n$, 
  it is possible to avoid the calculation of $4$ squares of lengths.
  The calculation formula for $C_n$ differs from simply-laced case by some parameter,
  which is fixed for all pairs $\{r, s\}$ with given extraspecial pair $\{r_1, s_1\}$.
\end{abstract}

\maketitle

\setcounter{tocdepth}{3}

\section{\bf Introduction} 
  In this paper, the new formulas for calculating structural constants
  for the root systems of types $B_n$, $C_n$ and $F_4$ are constructed.
  In these formulas several pairs of squares of lengths
  can be reduced. The time efficiency of these formulas for $B_n$ is the same, and for $C_n$
  is almost the same, as for simply-laced case.
  The main tools used in this paper are the concept of special and extraspecial pairs introduced
  by R.~Carter in \cite[\S4.2]{Ca72},
  and the inductive algorithm for calculating structural constants,
  see \cite[p.59-60]{Ca72}, \cite{CMT3}, \cite[$5^0$]{V04}.

  Let  $\{r, s\}$ be any special pair, and let $\{r_1, s_1\}$ be an extraspecial pair corresponding to $\{r, s\}$,
  see \S\ref{sec_extrasp}.
  We introduce the notion of a {\it quartet} - an ordered set of four roots  $\{r_1, r, s, s_1 \}$.
  The study of quartets has proven to be a useful tool for understanding the internal structure of root systems and
  for deriving formulas for calculating structural constants, see \S\ref{sec_main_res}.
  A detailed description of quartets for cases $B$, $C$ and $F$ are given in
  \S\ref{sec_Bn}-\S\ref{sec_F4}. As an example, Tables \ref{table_B6_quartets}-\ref{table_F4_quartets}
  provide lists of quartets for the cases $B_6$, $C_6$, $F_4$.

  The classification of quartets makes it possible to simplify the formulas for calculations.
  Quartet $q = \{r_1, r, s, s_1\}$ is said to be a {\it mono-quartet} if vectors $s-r_1$ and $r-r_1$
  are not roots at the same time.
  Quartet $q$ is said to be {\it simple} if it satisfies
  the following properties:
  if $s-r_1$ (resp. $r-r_1$) is a root, then $\abs{s-r_1} = \abs{s}$ (resp. $\abs{r-r_1} = \abs{r}$).
  By Lemma \ref{lem_Nab} formula for calculation the structure constant for the quartet $q$
  is as follows:
\begin{equation}
  \label{eq_lem_1_0}
  \footnotesize
\begin{split}
    N(r, & s)  =
    \frac{\abs{r+s}^2}{N(r_1, s_1)}
   \bigg ( \frac{N(s - r_1, r_1) N(s_1 - r, r) \abs{s_1 - r}^2}{\abs{s}^2\abs{s_1}^2} +
           \frac{N(r_1, r - r_1) N(s_1 - s, s) \abs{r - r_1}^2}{\abs{r}^2\abs{s_1}^2}  \bigg ).
\end{split}
\end{equation}
  If $q$ is a mono-quartet, then either $N(s-r_1, r_1) = 0$ or $N(r_1, r - r_1) = 0$,
  and in formula \eqref{eq_lem_1_0} one of the two terms in parentheses disappears.
  If $q$ is a simple quartet, then $2$ pairs of squared lengths in fraction \eqref{eq_lem_1_0}
  immediately cancel out,
  namely $\abs{s-r_1}^2/s^2 = 1$ and $\abs{r-r_1}^2/\abs{r}^2 = 1$.
  In this case, formula \eqref{eq_lem_1_0} is transformed to the following:
\begin{equation}
  \label{eq_Cn_final_0}
   N(r, s)  =
       \frac{\abs{r_1 + s_1}^2}{\abs{s_1}^2}
       \bigg ( \frac{N(s - r_1, r_1) N(s_1 - r, r) + N(r_1, r - r_1) N(s_1 - s, s)}{N(r_1, s_1)} \bigg ).
\end{equation}
Formula \eqref{eq_Cn_final_0} allows us to avoid calculating the squares of the lengths for four roots.
As it is shown in Theorems \ref{th_Bn_2} and \ref{th_Cn_2} any quartet in $B_n$ and $C_n$ is simple.
Then, for $B_n$ and $C_n$ we can use formula \eqref{eq_Cn_final_0}, and for $F_4$ we can use
formula \eqref{eq_lem_1_0}.  At last, by Theorem \ref{th_Bn_2},
for $B_n$ the roots $r_1 + s_1$ and $s_1$ have the same length,
and formula \eqref{eq_Cn_final_0} for $B_n$  is transformed as follows:
\begin{equation}
  \label{eq_Bn_final_0}
   N(r, s)  =
    \frac{N(s - r_1, r_1) N(s_1 - r, r) + N(r_1, r - r_1) N(s_1 - s, s)}{N(r_1, s_1)}.
\end{equation}
Formula \eqref{eq_Bn_final_0} can be applied to any simply-laced case, since the lengths of all roots
in formula \eqref{eq_lem_1_0} are reduced.
Note that in case $B_n$ there are roots of different lengths, but
formula \eqref{eq_Bn_final_0} does not depend on their lengths.

\subsection{The basic Carter formula}
   \label{sec_basic}
The formula underlying the calculation of structural constants
and from which formulas (1), (2) are derived is the Carter formula:
\begin{equation}
 \label{eq_intro_1}
 \begin{split}
   & \frac{N(\alpha, \beta)N(\gamma, \delta)}{\abs{\alpha+\beta}^2} +
     \frac{N(\beta, \gamma)N(\alpha, \delta)}{\abs{\beta+\gamma}^2} +
     \frac{N(\gamma, \alpha)N(\beta, \delta)}{\abs{\alpha+\gamma}^2} = 0,
  \end{split}
\end{equation}
if  \, $\alpha, \beta, \gamma, \delta \in \varPhi$, \, $\alpha + \beta + \gamma + \delta = 0$
and if no pair are opposite.
If one of sums in  \eqref{eq_intro_1} is not a root, the corresponding
structure constant is $0$, and the corresponding term in \eqref{eq_intro_1} disappears,
see \cite[Th. 4.1.2]{Ca72}.  For the simple-laced case, the lengths\footnotemark[1] of all roots are equal, and
eq. \eqref{eq_intro_1} looks significantly  simpler:
\footnotetext[1]{Note that throughout this article, the square of the length of any vectors is
the value of the quadratic Tits form on this vector, see \cite[Ch.~2]{St08}.}
\begin{equation}
 \label{eq_simple_case}
 \begin{split}
    N(\alpha, \beta)N(\gamma, \delta) + &  N(\beta, \gamma)N(\alpha, \delta) +  N(\gamma, \alpha)N(\beta, \delta) = 0, \\
    & \text{ if } \alpha + \beta + \gamma + \delta = 0.
  \end{split}
\end{equation}
\subsubsection{Vavilov's observation}
 \label{sec_Vavilov_obs}
In \cite{V04}, Vavilov observed that only two summands in \eqref{eq_simple_case} can be non-zero.
The question is, how things are in the multiply-laced case,
and what else can be said about the relations and formulas associated with structural constants
in the multiply-laced case? To answer this question, we used the so-called special and extraspecial pairs
of roots introduced by Carter, \cite[p.58]{Ca72}.

\subsection{Extraspecial pairs and quartets}
  \label{sec_extrasp}
The positive roots in any root system can be ordered in such a way that
a root with smaller height precedes a root of larger height and  the roots with the same height
are ordered in {\it lexicographic} rule. The symbol $\prec$ is used to denote
the relation ``precede". This ordering of positive roots is called {\it regular}, \cite[p.~67]{V04}.

An ordered pair $\{r, s\}$ of positive roots\footnotemark[2] is called a
{\it special pair} if $r + s \in \varPhi$ and  $0 \prec r \prec s$. An ordered pair $\{r_1, s_1\}$ of roots is called
{\it extraspecial} if $\{r_1, s_1\}$ is special and if for all special pairs $\{r, s\}$ with
$r + s = r_1 + s_1$ we have $r_1 \preceq r$. Extraspecial pairs are in one-ro-one correspondence
with the roots in  $\varPhi^+ - \Delta$, where $\varPhi^+$ is the set of the positive roots, $\Delta$
is the set of the simple roots, see \cite[p.~58-59]{Ca72}.
If $\{r, s\}$ is a pair of roots which is special but not extraspecial, then there is unique
extraspecial pair $\{r_1, s_1\}$ such that
\footnotetext[2]{Here, the special pairs are defined only for the positive roots.
For the transition from positive roots to arbitrary ones, see Remark \ref{rem_alg_posneg}.}
\begin{equation}
 \label{eq_extr_sp}
 r + s = r_1 + s_1.
\end{equation}
 The set of four roots $\{r,s, -r_1, -s_1\}$ satisfies \eqref{eq_intro_1}, and
\begin{equation}
  \label{eq_quadr}
 0 \prec  r_1 \prec r  \prec s \prec s_1.
\end{equation}

 The importance of extraspecial pairs is determined by the following proposition by Carter:
 \begin{proposition}\cite[Prop.~4.2.2]{Ca72}.
  \label{prop_estrasp}
  The signs of the structure constants $N(r,s)$ may be chosen arbitrarily for extraspecial pairs $\{r, s\}$,
  and then the structure constants for all pairs are uniquely determined.
\end{proposition}

 Let $\{r_1, s_1\}$ be an extraspecial pair, and $\{r, s\}$ be a special pair
 satisfying \eqref{eq_extr_sp} and \eqref{eq_quadr}.
 The ordered set of four roots  $\{r_1, r, s, s_1 \}$  is said to be a {\it Carter quartet} or just a {\it quartet}.
For examples of quartets and extraspecial pairs for $B_6$, $C_6$ and $F_4$,
see Tables \ref{table_B6_quartets}, \ref{table_C6_quartets} and \ref{table_F4_quartets}.

\subsection{The main results}
 \label{sec_main_res}
  Let $\{r_1, r, s, s_1 \}$ be a quartet in any root system.
  By \eqref{eq_extr_sp} the set of roots $\{r_1, -r, -s, s_1 \}$ satisfies the Carter formula \eqref{eq_intro_1}.
  By {\bf Lemma \ref{lem_Nab}}, formula \eqref{eq_intro_1}
  is transformed to formula \eqref{eq_lem_1_0}

\subsubsection{Types of quartets}

In what follows, the following definitions will be used:

\begin{itemize}
  \item A quartet  $\{r_1, r, s, s_1 \}$ is said to be a {\it mono-quartet}
  if vectors $s-r_1$ and $r-r_1$ are not roots at the same time.
  In this case, one of two summands in \eqref{eq_lem_1_0} disappears.
  \item A root system is said to be of the {\it Vavilov type} if any quartet is a mono-quartet,
   see \S\ref{sec_Vavilov_obs}.
  \item A quartet  $\{r_1, r, s, s_1 \}$ is said to be {\it simple quartet} if it
  satisfies the following properties:

    if $s - r_1$ is a root, then $\abs{s - r_1} = \abs{s}$; \;\;
    if $r - r_1$ is a root, then $\abs{r - r_1} = \abs{r}$.

\end{itemize}
 A simple quartet is not necessarily mono,
 such examples are found among $C_n$ quartets, see \S\ref{sec_Cn}.
 Conversely, a mono-quartet is not necessarily simple,
 examples of this are found among $F_4$ quartets, see \S\ref{sec_F4}.
 At last, all $B_n$ quartets are mono and simple, see \S\ref{sec_Bn}.

\subsubsection{Quartets for root system $B_n$}
  \label{sec_Bn_quartets}
In {\bf Theorem \ref{th_Bn_1}} it is proved that $B_n$ is of the Vavilov type.
In fact, this statement coincides with Vavilov's observation for the simply-laced case, see \S\ref{sec_Vavilov_obs}.
Further, in {\bf Theorem \ref{th_Bn_2}}, it is proved that
any quartet $\{r_1, r, s, s_1 \}$ in $B_n$ is simple,
and the roots $r_1 + s_1$ and $s_1$ have the same length.
This theorem allows us to simplify formula \eqref{eq_lem_1_0}:

\begin{corollary}
 \label{col_Bn_1}
In the root system $B_n$, formula \eqref{eq_lem_1_0} is simplified to formula \eqref{eq_Bn_final_0}.
In fact, relation \eqref{eq_Bn_final_0} is applicable
for calculating structure constants in the simply-laced case and in the case $B_n$.
\end{corollary}

\subsubsection{Quartets for root system $C_n$}
  \label{sec_Cn_quartets}
The properties of quartets for $B_n$ are not the same as for $C_n$,
but are close. By {\bf Theorem \ref{th_Cn_1}},

  (i) The root system $C_n$ is not of the Vavilov type.
  A quartet $\{r_1, r, s, s_1 \}$ is a mono-quartet if and only if $(r_1, s_1) \neq 0$.
  This condition does not depend on the pair $\{r, s\}$.

  (ii) The number of extraspecial pairs $\{r_1, s_1\}$ with  $(r_1, s_1) = 0$  is $n-1$,
  whereas the total number of extraspecial pairs is $n^2-n$.

  (iii) The inner product $(r_1, s_1) = 0$ if and only if
\begin{equation*}
  \abs{r_1 + s_1} = \sqrt{2}, \quad  \abs{r_1} = \abs{s_1} = 1, \;\text{ and }\; \frac{\abs{r_1 + s_1}^2}{\abs{s_1}^2} = 2.
\end{equation*}

By {\bf Theorem \ref{th_Cn_2}}, similarly to $B_n$, in the root system $C_n$, any quartet is simple.

\begin{corollary}
 \label{cor_Cn}
In the root system $C_n$, formula \eqref{eq_lem_1_0} is simplified to formula \eqref{eq_Cn_final_0}
\end{corollary}
Formula \eqref{eq_Cn_final_0} differs from the simply-laced case by the factor
$\varphi =\abs{r_1 + s_1}^2/\abs{s_1}^2$,
which depends only on the extraspecial pair $\{r_1, s_1\}$.
The parameter $\varphi$ takes values $1/2, \, 1, \, 2$.
For more details about parameter $\varphi$, see relation \eqref{eq_cases_fi}.

\subsubsection{Quartets for root system $F_4$}
  \label{sec_Fn_quartets}
The quartets in $F_4$ are not structured as well as in $B_n$ and $C_n$.
There are $48$ different quartets in $F_4$, see Table~\ref{table_F4_quartets}.
By {\bf Theorem \ref{th_F4_1}}, there are $38$ simple quartets in $F_4$: $30$ with $\abs{r_1} = \sqrt{2}$,
see Table~\ref{table_F4_quartets}, and $8$ with $\abs{r_1} = 1$,
see \eqref{eq_F4_8_quart}. These $38$ quartets the structure constants can be calculated by
formula \eqref{eq_Cn_final_0}.

\subsubsection{Tables}

The root systems $B_6$ and $C_6$ are used as an example, see \S\ref{sec_app}.
Both contain $36$ positive roots. The regular ordering of positive roots listed
in Table \ref{tab_B6C6}. Near each root (on the right) are the squares of the root lengths.
For a complete list of $B_6$ quartets (resp. $C_6$ quartets) containing $80$ quartets,
see  Table \ref{table_B6_quartets} (resp. Table \ref{table_C6_quartets}).
For a complete full list of $F_4$ quartets containing $48$ quartets, see  Table \ref{table_F4_quartets}.

\subsubsection{Implementation}
An implementation of several classes in Python for computing positive roots, 
extraspecial pairs, quartets, and structure constants can be found in \cite{Sk24}.

\subsection*{Acknowledgements}
I would like to thank A.~Elashvili and M.~Jibladze for useful discussions and valuable comments
in preparing this paper.

\section{\bf Structure constants in the multiply-laced case}

\subsection{Chevalley basis, Weyl basis and Chevalley theorem}
In \cite{Ch55} Chevalley constructed in a simple Lie algebra $\mathfrak{g}$
a basis with the property that all structure constants are integers.
Let
\begin{equation*}
 \label{eq_coroot}
    h_{\alpha} = \frac{2\alpha}{(\alpha, \alpha)}
\end{equation*}
be the {\it dual root} corresponding to $\alpha \in \varPhi$,
where $\varPhi$ is the root system corresponding to $\mathfrak{g}$,
$(\cdot, \cdot)$ is the corresponding Killing form\footnotemark[2].
Let $\mathfrak{h}$ be a Cartan subalgebra of
$\mathfrak{g}$. The set $\{h_{\alpha} | \alpha \in  \Delta \}$ forms the basis in $\mathfrak{h}$,
where $\Delta$ is a set of positive simple roots in $\varPhi$.
For the simple Lie algebra $\mathfrak{g}$ over $\mathbb{C}$, there exists a {\it Cartan decomposition}
\footnotetext[2]{The Killing form $(\cdot,\cdot)$ is defined on a Lie algebra.
When this form is considered on the root system,
it is also commonly called the Tits form, see \cite{K80}.}
\begin{equation*}
    \mathfrak{g} = \mathfrak{h} + \sum\limits_{\alpha \in \varPhi}\mathfrak{g}_{\alpha},
\end{equation*}
there $\{ \mathfrak{g}_{\alpha} \mid \alpha \in \varPhi \}$ are onedimensional root spaces, which satisfy
the following Lie relation:
\begin{equation*}
     [g_{\alpha}, g_{\beta}] = g_{\alpha+\beta},
\end{equation*}
For each root $\alpha$, choose a basis vector $e_{\alpha} \in \mathfrak{g}_{\alpha}$
such that
\begin{equation*}
     [e_{\alpha}, e_{\beta}] =
     \begin{cases}
        N(\alpha, \beta)e_{\alpha+\beta} \;\text{ if } \alpha+\beta \in \varPhi, \\
        0 \qquad\qquad\quad\; \text{ if } \alpha+\beta \not\in \varPhi.
     \end{cases}
\end{equation*}
where $N(\alpha, \beta)$ is some scalar.
Note that $N(\alpha, \beta)$ is assigned to $0$ if $\alpha + \beta$ is not a root.
The scalars $N(\alpha, \beta)$ determine the multiplication of elements
$e_{\alpha}$. The vectors $e_{\alpha}$ are called the {\it root vectors}, and scalars $N(\alpha, \beta)$
are called the {\it structure constant} of $\mathfrak{g}$. Note that $e_{\alpha}$ can be chosen up to coefficient,
and $N(\alpha, \beta)$ changes accordingly. The set
$\ \{h_{\alpha} \mid \alpha \in  \Delta;   e_{\alpha} \mid \alpha \in \varPhi \}$ forms the basis  in
$\mathfrak{g}$ which is called a {\it Weyl basis}.
The Weyl basis normalized in such a way that all the  structure constants are integers is called
the Chevalley basis.
Moreover, the following theorem holds:
\begin{theorem}[Chevalley, 1955]
 \label{th_Chevalley}
The root vectors $e_{\alpha}$ can be chosen so that
\begin{equation}
  \label{eq_rel_Nab}
  \begin{split}
     & [e_{\alpha}, e_{-\alpha}] = h_{\alpha} \text{ for all } \alpha \in \varPhi, \\
     & N(\alpha, \beta) = -N(-\alpha, -\beta)  \text{ for all } \alpha, \beta \in \varPhi.
  \end{split}
\end{equation}
 Let $\alpha, \beta \in \varPhi$ such that $\alpha + \beta \in \varPhi$ and $p$ be the greatest
integer such that $\beta-p\alpha \in \varPhi$.  Then one has
\begin{equation}
  \label{eq_pm_p_1}
   N(\alpha, \beta) = \pm(p+1),
\end{equation}
see \cite[Ch.~6,$\S$6]{Se01}, \cite[Th.~1]{Ch55}.
\end{theorem}

By eq. \eqref{eq_pm_p_1}
the structure constants $N(\alpha, \beta)$ take only the values $\pm{1}, \pm{2}, \pm{3}$. For the simply-laced
root systems, there cannot be a chain like this
\begin{equation}
  \label{eq_chain_3}
   \beta-\alpha, \;\; \beta, \;\; \beta + \alpha.
\end{equation}
Then, $p = 0$, and $N(a,b) = \pm{1}$. The chains like \eqref{eq_chain_3} exist for root systems
$F_4$, $B_n$ and $C_n$. Here, $p = 1$, and $N(a,b) = \pm{1}, \pm{2}$.  The chains of length $4$ exist
only for $G_2$, see \cite[Ch.6,~$\S$1,~$3^o$]{Bo02}. Then, $p = 2$ and $N(a,b)$ can take values $\pm{1}, \pm{2}, \pm{3}$.
There is just one catch left, how to determine the signs of the structure constants.

\subsection{Carter's theorem on structure constants}
The main properties of structure constants are listed in the following theorem \cite[Th.4.1.2]{Ca72}:
\begin{theorem}[Carter, 1972]
 \label{th_Carter}
   The structure constants of a simple Lie algebra $\mathfrak{g}$ over $\mathbb{C}$
satisfy the following relations:
\begin{equation}
 \label{eq_Nab_2}
   N(\alpha, \beta) = -N(\beta, \alpha), \qquad  \alpha, \beta \in \varPhi.
\end{equation}

\begin{equation}
  \label{eq_transf_pairs}
  \frac{N(\alpha, \beta)}{\abs{\gamma}^2} = \frac{N(\beta, \gamma)}{\abs{\alpha}^2} =
  \frac{N(\gamma, \alpha)}{\abs{\beta}^2}
\end{equation}
if \, $\alpha, \beta, \gamma \in \varPhi$ \, and \,  $\alpha + \beta + \gamma = 0$.

\begin{equation}
  \label{eq_3_memb}
   \frac{N(\alpha, \beta)N(\gamma, \delta)}{\abs{\alpha+\beta}^2} +
   \frac{N(\beta, \gamma)N(\alpha, \delta)}{\abs{\beta+\gamma}^2} +
   \frac{N(\gamma, \alpha)N(\beta, \delta)}{\abs{\alpha+\gamma}^2} = 0
\end{equation}
 if  \, $\alpha, \beta, \gamma, \delta \in \varPhi$, \, $\alpha + \beta + \gamma + \delta = 0$
 \, and if no pair are opposite. If one of sums in  \eqref{eq_3_memb} is not a root, the corresponding
 structure constant is 0, and the corresponding term in  \eqref{eq_3_memb} disappears.
\end{theorem}
One more relation we get by \eqref{eq_rel_Nab} and \eqref{eq_Nab_2}:
\begin{equation}
  \label{eq_neq2pos}
     N(-\alpha, -\beta) = N(\beta, \alpha).
\end{equation}

The relation \eqref{eq_transf_pairs} is due to J.~Tits \cite{T66}, see \cite{Cs18}.
Following B.~Casselman, \cite{Cs18}, the triple $\{\alpha, \beta, \gamma\}$  from \eqref{eq_transf_pairs}
is called the {\it Tits triple}.

\subsection{How to choose the signs of the structure constants?}

By using induction algorithm on the sum $r+s$ it can be see  that $N(r,s)$
is determined by the values of the structure constants on the extraspecial pairs, \cite[p.~60]{Ca72}.
See also \cite[$5^o$]{V04}, \cite{CMT3}. Let $\{r_1, r, s, s_1 \}$ be a quartet. Then
four roots $\{r_1, -r, -s, s_1 \}$ satisfy Theorem \ref{th_Carter} and \eqref{eq_3_memb}.
By \eqref{eq_3_memb} we have as follows:
\begin{equation*}
   \frac{N(r, s)N(-r_1, -s_1)}{\abs{r+s}^2} +
   \frac{N(s, -r_1)N(r, -s_1)}{\abs{s-r_1}^2} +
   \frac{N(-r_1, r)N(s, -s_1)}{\abs{r-r_1}^2} = 0
\end{equation*}

 By \eqref{eq_rel_Nab} we have $N(-r_1, -s_1) = -N(r_1, s_1)$ and by \eqref{eq_3_memb} we get the following expression
 for $N(r, s)$ :
\begin{equation}
  \label{eq_3_memb_2}
   N(r, s)  =
   \frac{\abs{r+s}^2}{N(r_1, s_1)}
   \bigg (\frac{N(s, -r_1)N(r, -s_1)}{\abs{s-r_1}^2} + \frac{N(-r_1, r)N(s, -s_1)}{\abs{r-r_1}^2} \bigg )
\end{equation}
 Note that $N(r_1, s_1) \neq 0$ because the pair $\{r_1, s_1\}$ is extraspecial, then
 the division by $N(r_1, s_1)$ in eq. \eqref{eq_3_memb_2} is acceptable.

 \begin{lemma}
  \label{lem_1}
    Let $\alpha, \beta \in \varPhi$ and $\alpha \prec \beta$.

   \rm{(i)} If $\alpha$ and $\beta$ are roots with the same height in $\varPhi$ then $\beta - \alpha$ is not a root.

   \rm{(ii)} If $\beta - \alpha$ is a root then $\beta - \alpha$ is positive.
 \end{lemma}
  \PerfProof (i) Since $\alpha \prec \beta$ then there is a coordinate $i$ such that $(\alpha)_i \neq (\beta)_i$.
  Assume that $(\alpha)_i > (\beta)_i$. Since $ht(\alpha) = ht(\beta)$ then there exists another coordinate $j$
  such that $(\alpha)_j < (\beta)_j$.  Thus, $\beta - \alpha$ has coordinates with different signs
  which cannot be.

  (ii) By (i) we have $ht(\beta) > ht(\alpha)$, then there exists coordinate $i$ such that $(\beta)_i > (\alpha)_i$,
  i.e,  $(\beta - \alpha)_i > 0$. Therefore, $\beta - \alpha > 0$. \qed
  \medskip

  Let us transform structure constants $N(s, -r_1)$, $N(r, -s_1)$, $N(-r_1, r)$, $ N(s, -s_1)$ from eq. \eqref{eq_3_memb_2}
  in such a way that calculations are performed only on positive roots. In addition, the obtained roots should have
  the height no greater than $r$ and $s$.

\begin{lemma}
 \label{lem_Nab}
  Let $\{r_1, r, s, s_1 \}$ be a quartet in any multiply-laced root system.
  By \eqref{eq_extr_sp}, the roots $\{r_1, r, s, s_1 \}$ satisfy the Carter formula \eqref{eq_3_memb},
  which can be transformed as follows:
\begin{equation}
  \label{eq_lem_1}
\begin{split}
    N(r, & s)  =
    \frac{\abs{r+s}^2}{N(r_1, s_1)}
   \bigg ( \frac{N(s - r_1, r_1) N(s_1 - r, r) \abs{s_1 - r}^2}{\abs{s}^2\abs{s_1}^2} +
           \frac{N(r_1, r - r_1) N(s_1 - s, s) \abs{r - r_1}^2}{\abs{r}^2\abs{s_1}^2}  \bigg )
\end{split}
\end{equation}
\end{lemma}

\PerfProof  First, the following relations hold:
  \begin{equation}
    \label{eq_Nab_transf}
    \begin{split}
      & N(s, -r_1) =  N(s - r_1, r_1) \frac{\abs{s - r_1}^2}{\abs{s}^2}, \qquad
        N(r, -s_1) =  N(s_1 - r, r) \frac{\abs{s - r_1}^2}{\abs{s_1}^2}, \\
      & N(-r_1, r) =  N(r_1, r - r_1)\frac{\abs{r-r_1}^2}{\abs{r}^2}. \qquad
        N(s, -s_1) =  N(s_1 - s, s) \frac{\abs{r - r_1}^2}{\abs{s_1}^2}.
    \end{split}
  \end{equation}

  Let us prove the first relation of \eqref{eq_Nab_transf}.  By \eqref{eq_neq2pos} $N(s, -r_1) = N(r_1, -s)$.
  Since the triple of roots $\{r_1, -s, s - r_1\}$ forms the Tits triple satisfying eq. \eqref{eq_transf_pairs},
  we have
  \begin{equation*}
     \frac{N(s, -r_1)}{\abs{s - r_1}^2} = \frac{N(r_1, -s)}{\abs{s - r_1}^2} =
       \frac{N(s - r_1, r_1)}{\abs{s}^2}.
  \end{equation*}
  Since $0 \prec s - r_1$, by Lemma \ref{lem_1} we have as follows:
  if $s - r_1$ is a root then  $s - r_1$ is positive.
  Otherwise, $N(s - r_1, r_1) = 0$.

  The remaining three relations are proved similarly.
  Note that in the second (resp. forth) relation of \eqref{eq_Nab_transf}, we use also the relation
  $s_1 - r = s - r_1$ (resp. $s - s_1 = r_1 - r$).
  Substituting relations \eqref{eq_Nab_transf} into \eqref{eq_3_memb_2} we obtain the
  required formula \eqref{eq_lem_1} for $N(r,s)$. \qed
~\\

\begin{remark}
  \label{rem_alg_posneg}
{\rm
Consider calculating structure constants for the pairs $\{\alpha, \beta \}$ such that $\alpha > 0$,  $\beta < 0$.
If $\alpha + \beta > 0$, then by \eqref{eq_transf_pairs} and \eqref{eq_neq2pos}, we have
\begin{equation}
 \label{eq_negpos_2}
  \frac{N(\alpha, \beta)}{\abs{\alpha + \beta}^2} \; = \; \frac{N(\beta, -\alpha-\beta)}{\abs{\alpha}^2} \; = \;
  \frac{N(\alpha + \beta, -\beta)}{\abs{\alpha}^2}.
\end{equation}
 In \eqref{eq_negpos_2}, the last expression $N(\alpha + \beta, -\beta)$ already has two positive roots
 as input parameters. If $\alpha + \beta < 0$, then by \eqref{eq_transf_pairs}
\begin{equation}
 \label{eq_negpos_3}
  \frac{N(\alpha, \beta)}{\abs{\alpha + \beta}^2} \; = \; \frac{N(-\alpha-\beta, \alpha)}{\abs{\beta}^2}.
\end{equation}
In \eqref{eq_negpos_3}, the last expression $N(-\alpha - \beta, \alpha)$ again has two positive roots
as input parameters.
}
\end{remark}
\section{\bf Roots and extraspecial pairs in the multiply-laced case}
\subsection{The positive roots in $B_n$ and $C_n$}
Each of the root systems $B_n$ and $C_n$ contains $n^2$ positive roots.
In $B_n$ there are $n$ roots of length $1$ and $n^2 - n$ roots of length $\sqrt{2}$,
and, conversely, in $C_n$ there are $n$ roots of length $\sqrt{2}$ and $n^2 - n$ roots of length $1$.
\subsubsection{Root system $B_n$}
\begin{equation*}
 \begin{split}
   & \text{ Positive roots: } 
   \varepsilon_i \; (1 \le i \le n), \quad \varepsilon_i \pm  \varepsilon_j \; (1 \le i < j \le n). 
    \qquad\qquad\qquad \qquad\qquad\qquad \qquad\qquad\qquad \\
   & \text{ Basis: } 
     \{ \alpha_1 = \varepsilon_1 - \varepsilon_2, \;\; \alpha_2 = \varepsilon_2 - \varepsilon_3, \;\; \dots, \;\;
    \alpha_{n-1} = \varepsilon_{n-1} - \varepsilon_n, \;\; \alpha_n = \varepsilon_n \}, \\
 \end{split} 
\end{equation*}
see \cite[Plate II]{Bo02}.
Let $\beta_{ij} := \varepsilon_i - \varepsilon_{j+1}$, $\alpha_{ij} := \varepsilon_i + \varepsilon_{j}$.
See Table \ref{tab_Bn} for the positive roots of $B_n$. 

\begin{table}[!ht]
\centering
\renewcommand{\arraystretch}{1.2}
  \begin{tabular} {|c|c|c|}
\hline  
  \footnotesize
  $\begin{array}{l}
   \setlength\arraycolsep{1pt}
   \begin{matrix}
     \varepsilon_i =  & 0 & \dots & 0 & 1 & \dots & 1 \\
                &   &       &   & i   \\
    \end{matrix}  \\
  \end{array}$ & 
  \footnotesize
  $\begin{array}{l}
  \setlength\arraycolsep{1pt}
   \begin{matrix}
    \beta_{ij} =  & 0 & \dots & 0 & 1 & \dots & 1 & 0 & \dots & 0 \\
                  &   &       &   & i &       & j      \\
  \end{matrix} \\
   \end{array}$ &
   \footnotesize
   $\begin{array}{l}
  \setlength\arraycolsep{1pt}
   \begin{matrix}
    \alpha_{ij} =  & 0 & \dots & 0 & 1 & \dots & 1 & 2 & \dots & 2 \\
                   &   &       &   & i &       &   & j   \\
  \end{matrix} \\
   \end{array}$ \\   
 \hline    
   \footnotesize  $1 \leq i \leq n$  & \footnotesize $1 \leq i \leq j \leq n-1$  \footnotesize & $1 \leq i < j \leq n$ \\
 \hline
    \footnotesize $n$ roots, \; $\abs{\varepsilon_i} = 1$ &
    \footnotesize $ \displaystyle\frac{n(n-1)}{2}$ roots, \; $\abs{\beta_{ij}} = \sqrt{2}$ &
    \footnotesize $ \displaystyle\frac{n(n-1)}{2}$  roots, \; $\abs{\alpha_{ij}} = \sqrt{2}$ \\
 \hline  
  \end{tabular}
  ~\\ \vspace{2mm}
  \caption{The positive roots of $B_n$ in the basis $\{\alpha_1, \dots, \alpha_n \}$}
  \label{tab_Bn}
\end{table}

\subsubsection{Root system $C_n$}
\begin{equation*}
 \begin{split}
   & \text{ Positive roots: }
   2\varepsilon_i \; (1 \le i \le n), \quad \varepsilon_i \pm  \varepsilon_j \; (1 \le i < j \le n).
    \qquad\qquad\qquad \qquad\qquad\qquad \qquad\qquad\qquad \\
   & \text{ Basis: }
     \{ \alpha_1 = \varepsilon_1 - \varepsilon_2, \;\; \alpha_2 = \varepsilon_2 - \varepsilon_3, \;\; \dots, \;\;
    \alpha_{n-1} = \varepsilon_{n-1} - \varepsilon_n, \;\; \alpha_n = 2\varepsilon_n \}, \\
 \end{split}
\end{equation*}
see \cite[Plate III]{Bo02}. 
Let $\delta_{ij} := \varepsilon_i - \varepsilon_{j+1}$, \; $\gamma_{ij} := \varepsilon_i + \varepsilon_{j}$.  

\begin{table}[!ht]
\centering
\renewcommand{\arraystretch}{1.2}
  \begin{tabular} {|c|c|c|}
\hline
  \footnotesize
  $\begin{array}{l}
   \setlength\arraycolsep{1pt}
   \begin{matrix}
    2\varepsilon_i =   & 0 & \dots & 0 &  2  & \dots & 2  & 1 \\
                &   &       &   &  i  &           \\
  \end{matrix}  \\
  \end{array}$ &
  \footnotesize
  $\begin{array}{l}
  \setlength\arraycolsep{1pt}
   \begin{matrix}
    \delta_{ij} =  & 0 & \dots & 0 & 1 & \dots & 1 & 0 & \dots & 0  \\
                   &   &       &   & i &       & j &      \\
  \end{matrix} \\
   \end{array}$ &
   \footnotesize
   $\begin{array}{l}
  \setlength\arraycolsep{1pt}
   \begin{matrix}
    \gamma_{ij} =  0 & \dots & 0 & 1 & \dots & 1 & 2 & \dots & 2 & 1 \\
                     &       &   & i &       &   & j &    \\
  \end{matrix} \\
   \end{array}$ \\
 \hline
   \footnotesize  $1 \leq i \leq n$  & \footnotesize $1 \leq i \leq j \leq n-1$  \footnotesize & $1 \leq i < j \leq n$ \\
 \hline
    \footnotesize $n$ roots, \; $\abs{2\varepsilon_i} = \sqrt{2}$ &
    \footnotesize $ \displaystyle\frac{n(n-1)}{2}$ roots, \; $\abs{\delta_{ij}} = 1$ &
    \footnotesize $ \displaystyle\frac{n(n-1)}{2}$  roots, \; $\abs{\gamma_{ij}} = 1$ \\
 \hline
  \end{tabular}
  ~\\ \vspace{2mm}
  \caption{The positive roots of $C_n$ in the basis $\{\alpha_1, \dots, \alpha_n \}$}
  \label{tab_Cn}
\end{table}

Note that in Table \ref{tab_Cn}, $\{2\varepsilon_n = 0 \dots 0 1\}$, and $\gamma_{i,n} = \{0 \dots 0 1\dots 1\}.$
\medskip

\begin{lemma}[H72, Ch.10.2]
  \label{lem_Humphri}
  If $\alpha$ is a positive but non-simple root in $\varPhi$, then $\alpha - \beta$
  is a positive root for some simple root $\beta$. \qed
\end{lemma}

\begin{corollary}
  \label{col_extrasp}
  {\rm(i)} For any positive non-simple root $\alpha \in \varPhi$, there exist unique simple root $\beta$
  such that $\{\beta, \alpha - \beta\}$ is the extraspecial pair.

  {\rm(ii) }Each of the root systems $B_n$ and $C_n$ has $n^2 - n$ extraspecial pairs.
\end{corollary}
  \PerfProof (i) If both $\beta$ and $\alpha - \beta$ are simple then such a pair is unique for given $\alpha$.
  Now, let $\alpha - \beta$ be a non-simple root.
  The simple root $\beta$ precedes to $\alpha - \beta$, since $\Ht(\beta) < \Ht(\alpha - \beta)$.
  Let $\{\beta_1, \alpha - \beta_1 \}$ and $\{\beta_2, \alpha - \beta_2 \}$ be two special pairs,
  where $\beta_1$ and $\beta_2$ are simple.
  If $i$ (resp. $j$) is the non-zero coordinate of $\beta_1$ (resp. $\beta_2$) and $i < j$
  then $\beta_1 \prec \beta_2$. In addition, $i$th coordinate is the first where
  $(\alpha - \beta_1)_i < (\alpha - \beta_2)_i$, i.e., $\alpha - \beta_2 \prec \alpha - \beta_1$.
  Thus, $\beta_1 \prec \beta_2 \prec  \alpha - \beta_2 \prec \alpha - \beta_1$.
  Here, $\{\beta_1,  \alpha - \beta_1\}$  is extraspecial. \qed

  (ii) There are $n^2$ positive roots in each of the root systems $B_n$  and $C_n$.
  By (i) each of $B_n$ and $C_n$ contains $n^2 - n$ extraspecial pairs,
  and for each extraspecial pair the first root is simple. \qed

\subsection{How to get extraspecial pairs?}
  \label{sec_alg_extrasp}
By Corollary \ref{col_extrasp} the number of extraspecial pairs is equal to the number of
non-simple roots. Build a dictionary, where for any  extraspecial pair $\{r_1, s_1\}$,
an entry with a key equal to the root $r_1 + s_1$ will be made.
The entry corresponding $key = \gamma$ is the list of all special pairs $\{r, s\}$
such that $r+s =\gamma$.

Assigning structure constants to extraspecial pairs will be made using the Chevalley relation \eqref{eq_pm_p_1}.
If $s_1 - r_1$ is not a root, we set $N(s_1, r_1) := 1$. Else, if
$s_1 - r_1$ is a root,  but $s_1 - 2r_1$ is not a root, we set $N(s_1, r_1) := 2$.
At last, if $s_1 - r_1$ is a root and $s_1 - 2r_1$ is a root, we assign $N(s_1, r_1) := 3$.

See examples of extraspecial pairs and quartets  for $B_6$, $C_6$ and $F_4$
in Tables \ref{table_B6_quartets}, \ref{table_C6_quartets} and \ref{table_F4_quartets}.

\subsection{Extraspecial pairs and quartets in $B_n$}
   \label{sec_Bn}
 Let $\{r_1, s_1\}$ be an extraspecial pair in $B_n$.
For a simple root $r_1$  there are two possible cases:
\begin{equation*}
 \begin{split}
   (1) & \;\; r_1 \text{ is one of } \beta_{i,i}, \; i = 1, \dots, n-1. \text{ Then } \abs{r_1} = \sqrt{2}.  \\
   (2) & \;\; r_1 \text{ is } \varepsilon_n. \text{ Then } \abs{r_1} = 1.
 \end{split}
\end{equation*}

\begin{proposition}[Case $B_n$, $\abs{r_1} = \sqrt{2}$]
  \label{prop_Bn_1}
  Consider the root system $B_n$.

  {\rm(i)} For any extraspecial pair $\{r_1,  s_1\}$  with $\abs{r_1} = \sqrt{2}$, we have $(r_1,  s_1) = -1$.

  {\rm(ii)}  Let $q = \{r_1, r, s, s_1 \}$ be a quartet with $\abs{r_1} = \sqrt{2}$.
    Then, $q$ is a mono-quartet.
\end{proposition}

  \PerfProof (i) In the root system $B_n$, the inner product of two non-orthogonal roots
  can be only $\pm{1}$.  For $r_1 \in \{ \beta_{1,1},\dots, \beta_{n-1,n-1} \} $,
  the length of $\abs{r_1}$ is $\sqrt{2}$. Then,
  \begin{equation*}
        (r_1,  s_1)  =  (r_1, r_1 + s_1) - 2 < 0, \text{ i.e. } (r_1, s_1) = -1.
  \end{equation*}

  (ii) Suppose that both $s-r_1$ and $r-r_1$ are roots. Let us show, that in this case
  $(s, r_1) = \pm{1}$, and  $(r, r_1) = \pm{1}$.
  Really, if $s \perp r_1$, then
  \begin{equation*}
    \label{eq_orthog_1}
    \abs{s-r_1}^2 = \abs{s}^2 + \abs{r_1}^2 = \abs{s}^2 + 2,
  \end{equation*}
  which contradicts relation  $\abs{s-r_1}^2 \leq 2$. So, $(s, r_1) = \pm{1}$,
  and, similarly, $(r, r_1) = \pm{1}$. The inner product $(r + s, r_1)$ takes the values $0$ or $\pm{2}$.
  Because $r_1 + s_1 = r+s$, we have $(r_1 + s_1, r_1) = 2 + (r_1, s_1)$. Thus,
  $(r_1, s_1)$ is $0$ or $-2$, which contradicts (i). \qed

  \begin{proposition}[Case $B_n$, $r_1 = \varepsilon_n$]
    \label{prop_Bn_2}
    In the root system $B_n$,
    there are two extraspecial pairs $\{r_1, s_1\}$
    with $r_1 = \varepsilon_n$. In both cases, there is no
    other special pair $\{r, s \}$ such that $r + s = r_1 + s_1$,
    In other words, there is no quartet with  $r_1  = \varepsilon_n$.
  \end{proposition}
  \PerfProof
 Let us find the roots $s_1$ that, together with $r_1$, can form an extraspecial pair $\{r_1, s_1\}$.
 There are only two possible candidates for extraspecial pairs containing $\varepsilon_n$:

 (1)  Let $\{r_1, s_1\} = \{\varepsilon_n, \varepsilon_i\}$, where $i < n$.  For  $1  \leq i \leq n-2$,
 \begin{equation*}
 \begin{array}{l}
  \setlength\arraycolsep{1pt}
  \begin{matrix}
   r_1 + s_1 =  \varepsilon_n + \varepsilon_i =  & \{0 & \dots & 0 & 1 & \} + \\
                                     &
  \end{matrix}
  \begin{matrix}
     & \{0 & \dots & 0 & 1 & \dots & 1\} = \\
     &     &       &   & i
  \end{matrix}
  \begin{matrix}
     & \{0 & \dots & 0 & 1 & \dots & 1 & 1 & 2\}.\\
     &     &       &   & i
  \end{matrix}
 \end{array}
\end{equation*}
In this case, there is another expansion:
\begin{equation}
 \label{eq_another_dec_1}
\begin{array}{l}
\setlength\arraycolsep{1pt}
 \begin{matrix}
   r_1 + s_1 = \beta_{i,i} + \alpha_{i+1,n} = \{0  \dots  0 &  1 & 0  \dots   0\} +   \\
                                                            &  i &                        \\
  \end{matrix}
\begin{matrix}
    \{0 \dots & 0 & 1 & 1 \dots  1  2\}. \\
              & i &   \\
  \end{matrix}
 \end{array}
\end{equation}
So, $\beta_{i,i} \prec \varepsilon_n \prec \varepsilon_i \prec \alpha_{i+1,n}$.
Therefore, the pair $\{r_1, s_1\} = \{\varepsilon_n, \varepsilon_i \}$ is not
extraspecial.

For  $i = n-1$, we have the following  pair
$\{\varepsilon_n, \varepsilon_{n-1}\}$, which is extraspecial:
 \begin{equation}
   \label{eq_extrsp_Bn_1}
 \begin{array}{l}
  \setlength\arraycolsep{1pt}
  \begin{matrix}
   r_1 + s_1 = \varepsilon_n + \varepsilon_{n-1} = \{0  \dots  0 & 1 \} +  & \{0 & \dots & 0 & 1 & 1\}.\\
  \end{matrix}
 \end{array}
\end{equation}
However, there is no another pair $\{r,s\}$ such that $r+s = \varepsilon_n + \varepsilon_{n-1}$.
~\\

(2) Let $\{r_1, s_1\} = \{\varepsilon_n, \beta_{i,n-1}\}$, where $i < n$.
By \eqref{tab_Bn} we have
\begin{equation*}
 \begin{array}{l}
  \setlength\arraycolsep{1pt}
  \begin{matrix}
   r_1 + s_1 = \varepsilon_n + \beta_{i,n-1} = \{0  \dots  0 & 1 \} + & \{0 & \dots 0 & 1 & \dots & 1 & 0\}.\\
                                                       &        &     &         & i  \\
  \end{matrix}
 \end{array}
\end{equation*}
As above in \eqref{eq_another_dec_1}, for $i = 1,\dots, n-2$,
the root $\varepsilon_n + \varepsilon_i$ has another expansion:
\begin{equation*}
 \begin{array}{l}
  \setlength\arraycolsep{1pt}
  \begin{matrix}
    r_1 + s_1 = \beta_{i,i} + \varepsilon_{i+1} = \{0  \dots  0 &  1 & 0  \dots   0\} +  \{0 \dots & 0 & 1   &  \dots    1\}. \\
                                                       &  i &                             & i &   \\
  \end{matrix}
 \end{array}
\end{equation*}
Here, $\beta_{i,i} \prec \varepsilon_n \prec \beta_{i, n-1} \prec \varepsilon_{i+1}$.
Again, the pair $\{\varepsilon_n, \beta_{i,n-1} \}$ is not extraspecial.

For $s_1 = \beta_{n-1,n-1}$, we get an extraspecial pair $\{\beta_{n-1,n-1}, \varepsilon_n \}$:
 \begin{equation}
   \label{eq_extrsp_Bn_2}
 \begin{array}{l}
  \setlength\arraycolsep{1pt}
  \begin{matrix}
   r_1 + s_1 =  \beta_{n-1,n-1} + \varepsilon_n = \{0  \dots  0 & 1 & 0 \} +  & \{0 & \dots & 0 & 1\}.\\
  \end{matrix}
 \end{array}
\end{equation}
However, there is no another pair $\{r,s \}$ such that $r+s = \beta_{n-1,n-1} + \varepsilon_n$.
Thus, there are only two extraspecial pairs containing $\varepsilon_n$, see \eqref{eq_extrsp_Bn_1}
and \eqref{eq_extrsp_Bn_2}.
\qed
~\\

To summarize, we get the following
\begin{theorem}
  \label{th_Bn_1}
In the root system $B_n$, any quartet is a mono-quartet.
Thus, $B_n$ is of the Vavilov type.
\end{theorem}

It turns out that we can say something more precise about the lengths of $s - r_1$ and $r - r_1$.
\begin{theorem}
  \label{th_Bn_2}
  In the root system $B_n$, any quartet $\{r_1, r, s, s_1 \}$ is simple, and
  the roots $r_1 + s_1$ and $s_1$ have the same length.
\end{theorem}
\PerfProof
  Let $\abs{r_1} = 1$. Since $r_1$ is the simple root, by \eqref{tab_Bn} $r_1 = \varepsilon_n$.
By Proposition \ref{prop_Bn_2}, there is no quartet with $r_1 = \varepsilon_n$.

  Let $\abs{r_1}^2 = 2$, then $\abs{s - r_1}^2 = \abs{s}^2 + 2 - 2(s,r_1)$.
Since $s - r_1$ is a root, we have $(s,r_1) = 0, \pm{1}$.
If $(s,r_1) \leq 0$, then $\abs{s - r_1}^2 > 2$, which cannot happen.
Thus, $(s,r_1) = 1$, and $\abs{s - r_1}^2 = \abs{s}^2$. Similarly, if $r - r_1$ is a root,
then $\abs{r - r_1}^2 = \abs{r}^2$.

  The simple root $r_1$ can be either $\beta_{i,i}$, where $i < n$ or $\varepsilon_n$,
By Proposition \ref{prop_Bn_2}, there is no quartet with $r_1 = \varepsilon_n$.
For $r_1 = \beta_{i,i}$, we have $\abs{r_1}^2 = 2$
and $(r_1, s_1) = -1$, see Proposition \ref{prop_Bn_1}.
Then, $\abs{r_1+ s_1}^2 = \abs{r_1}^2 - 2 + \abs{s_1}^2 = \abs{s_1}^2$. \qed

\subsection{Extraspecial pairs and quartets in $C_n$}
  \label{sec_Cn}

Let $\{r_1, s_1 \}$ be an extraspecial pair in $C_n$.
In the case of $C_n$, dual to the root system of $B_n$,
for a simple root $r_1$ two cases are possible:
\begin{equation*}
 \begin{split}
   (1) & \;\; r_1 \text{ is one of } \delta_{i,i}, \; i = 1, \dots, n-1. \text{ Then } \abs{r_1} = 1.  \\
   (2) & \;\; r_1 \text{ is } 2\varepsilon_n. \text{ Then } \abs{r_1} = \sqrt{2}.
 \end{split}
\end{equation*}

\begin{proposition}[Case $C_n$, $\abs{s_1} = \sqrt{2}$ or $\abs{r_1} = \sqrt{2}$.]
  \label{prop_Cn_1}
Consider the root system $C_n$.

{\rm(i)}
 For each $\varepsilon_i$ with $i < n$, there exists the only extraspecial pair
 $\{r_1, s_1\}$ such that $s_1 = 2\varepsilon_i$, namely   $\{r_1, s_1\} = \{ \delta_{i-1,i-1}, 2\varepsilon_i\}$.
 In this case any quartet is a mono-quartet.

{\rm(ii)} There is no quartet $\{r_1, r, s, s_1 \}$ with $s_1 = 2\varepsilon_n$.

{\rm(iii)} There is no quartet $\{r_1, r, s, s_1 \}$ with $r_1 = 2\varepsilon_n$.

\end{proposition}

\PerfProof
(i) The root $s_1$ can be of length $\sqrt{2}$ only for $s_i = 2\varepsilon_i$, $i = 1,\dots, n$.
Consider case $s_1 = 2\varepsilon_i$, $i < n$.
Since $r_1$ is a simple root, there is only $r_1 = \delta_{i-1,i-1}$,
which forms an extensional pair. Here, $r_1 + s_1 = \gamma_{i-1,i}$, see Table \ref{tab_Cn}:
\begin{equation}
  \label{eq_Cn_prop_1}
  \footnotesize
  \setlength\arraycolsep{1pt}
  \begin{matrix}
     r_1 = \delta_{i-i,i-1} =  & \{0 & \dots & 0 & 1 &  0  & \dots & 0\} \\
                        &   &       &   &   &  i         \\
  \end{matrix}
  \qquad
  \begin{matrix}
     s_1 = 2\varepsilon_i =  & \{0 & \dots & 0 &  2  & \dots & 2  & 1 \}\\
                &   &       &   &  i  &           \\
  \end{matrix}  \\
  \qquad
  \begin{matrix}
     \gamma_{i-1,i} =  & \{0 & \dots & 1 &  2  & \dots & 2  & 1 \}\\
                &   &       &   &  i  &           \\
  \end{matrix}  \\
\end{equation}
 We have $\abs{r_1}$ = 1, $\abs{s_1} = \sqrt{2}$, and $\abs{r_1 + s_1} = 1$. The inner product $(r_1, s_1)$ is $-1$.
 Another expansion $r+s = \gamma_{i-1, i}$ can only be such that the coordinate
 $i-1$ of one of the pair $\{r, s\}$ is equal to $1$, and the other $0$.
 If $(s)_{i-1} = 0$, then $(s - r_1)_{i-1}$ is $-1$.
 Since $r_1 \prec s$, the vector $s - r_1$ has also positive coordinates,
 so $s - r_1$ is not a root. Similarly, if $(r)_{i-1} = 0$,
 $r - r_1$ is not a root.

(ii) Because $r_1$ is a simple root and $s_1 = 2\varepsilon_n$
then the only possibility for $r_1$ may be $\delta_{n-1, n-1}$, see \eqref{eq_Cn_prop_1}.
However, the root $\gamma = 2\varepsilon_n + \delta_{n-1, n-1}$ has no other expansion into the sum of
positive roots:
\begin{equation*}
 \footnotesize
  r_1 = \delta_{n-1, n-1} = \{0 \dots 1 0 \}, \qquad s_1 = 2\varepsilon_n = \{0 \dots 0 1\}, \qquad
      \gamma = \{0 \dots 1 1 \}.
\end{equation*}

(iii) If $r_1 = 2\varepsilon_n$, then $s_1$ is $\delta_{i,n-1}$ and $\gamma =  2\varepsilon_n + \delta_{i,n-1}$:
\begin{equation*}
 \footnotesize
  \setlength\arraycolsep{1pt}
  \begin{matrix}
   r_1 =  2\varepsilon_n =  & \{0  &  \dots & 0  & 1\}, \\
                      &    &        &    &  \\
  \end{matrix}  \qquad
  \begin{matrix}
   s_1 =  \delta_{i,n-1} =  & \{0 &  \dots & 0 & 1 &  \dots  & 1 & 0 \}, \\
                      &   &        &   & i &         &   &  \\
  \end{matrix}  \qquad
  \begin{matrix}
     \gamma =  & \{0 & \dots & 0 &  1  & \dots   & 1 \}. \\
               &   &       &   &  i  &         \\
  \end{matrix}  \\
\end{equation*}
Another expansion of $\gamma$ is as follows:
 \begin{equation*}
 \footnotesize
 \begin{array}{l}
  \setlength\arraycolsep{1pt}
  \begin{matrix}
   \gamma = \delta_{i,i} + \gamma_{i+1,n} =  &  \{  0 \dots  0 & 1 & 0 \dots 0 \} \; + \;
           & \{ 0  \dots & 0 &  1 &  \dots  1 \}. \\
                                                  &                 & i &
           &             & i  &   &    \\
  \end{matrix}
 \end{array}
\end{equation*}
Thus, the pair $\{r_1, s_1\} = \{2\varepsilon_n, \delta_{i,n-1} \}$  is not extraspecial.
\qed
\medskip

\begin{corollary}[Case $C_n$]
 \label{col_Cn_1}
 Let $\{r_1, s_1\}$ be the extraspecial root in $C_n$
 with $\abs{r_1} = \sqrt{2}$ or $\abs{s_1} = \sqrt{2}$. 
 Then, quartets exist only for
 $\{r_1, s_1\} = \{\delta_{i-1,i-1}, 2\varepsilon_i \}$.

{\rm(i)} In this case, the following relations take place:
  \begin{equation}
     \abs{r_1 + s_1} = 1, \quad  \abs{r_1} = 1, \quad  \abs{s_1} = \sqrt{2}, \quad (r_1, s_1) = -1 ~.
  \end{equation}

{\rm(ii)} Only two special pairs $\{r, s\}$ joined together with $\{r_1, s_1\}$ can form a quartet
 :
  \begin{equation}
   \label{en_decomp_rs_2}
   \begin{split}
     & \{r,s\} = \{\delta_{i-1,k}, \gamma_{i,k}\}, \quad \gamma_{i-1,i} = \delta_{i-1,k} + \gamma_{i,k+1},  \\
     & \{r,s\} = \{\delta_{i,k}, \gamma_{i-1,k}\}, \quad \gamma_{i-1,i} = \delta_{i,k} + \gamma_{i-1,k+1}.
   \end{split}
  \end{equation}
  Here, $\abs{r} = \abs{s} = 1$.
\end{corollary}
\PerfProof
 (i) Since $\gamma_{i-1,i-1} = r_1 + s_1$, the length of $r_1 + s_1$ is $1$, see  \eqref{eq_Cn_prop_1}
 and  Table \ref{tab_Cn}.  Then,
 \begin{equation*}
   \abs{r_1 + s_1}^2 = 3 + 2(r_1, s_1) = 1, \text{ so } (r_1, s_1) = -1.
 \end{equation*}

 (ii) The pairs $\{r,s\}$ in expansion \eqref{en_decomp_rs_2} are as follows:
 \begin{equation*}
  \label{eq_Cn_col_1}
  \footnotesize
  \setlength\arraycolsep{1pt}
 \begin{split}
  & \begin{matrix}
    r = \delta_{i-i,k} =   &\{0 & \dots & 0 & 1 &  1  & \dots & 1  & 0 \dots 0\}, \\
                            &    &       &   &   &  i  &       & k       \\
  \end{matrix}
   \quad
  \begin{matrix}
     s = \gamma_{i,k+1} =   &\{0 & \dots &  0 &  1  & \dots & 1  & 2 \dots 2 1\}, \; \text{ or} \\
                            &    &       &    &  i  &       & k       \\
  \end{matrix}  \\
 & \begin{matrix}
    r = \delta_{i,k}   =   &\{0 & \dots & 0 & 0 &  1  & \dots & 1  & 0 \dots 0\}, \\
                           &    &       &   &   &  i  &       & k       \\
  \end{matrix}
   \quad
  \begin{matrix}
     s = \gamma_{i-i,k+1} =   &\{0 & \dots &  1 &  1  & \dots & 1  & 2 \dots 2 1\}. \\
                            &    &       &    &  i  &       & k       \\
  \end{matrix}
\end{split}
\end{equation*}
In all cases, by Table \ref{tab_Cn} $\abs{r} = \abs{s} = 1$. \qed
\medskip

In Proposition \ref{prop_Cn_1} and Corollary \ref{col_Cn_1},
we considered the case where the length at least one of the pair $\{r_1,  s_1\}$
differs from $1$. In what follows, we assume, on the contrary,  $\abs{r_1} = \abs{s_1} = 1$.

\begin{proposition}[Case $C_n$, $\abs{r_1} = \abs{s_1} = 1$; $\abs{r} = \sqrt{2}$ or $\abs{s} = \sqrt{2}$]
  \label{prop_Cn_2}
 Let  $q = \{r_1, r, s, s_1\}$ be a quartet in $C_n$ such that $\abs{r_1} = \abs{s_1} = 1$
and $\abs{r} = \sqrt{2}$ (resp. $\abs{s} = \sqrt{2}$). Then $q$ is a mono-quartet.
\end{proposition}
 \PerfProof
 Assume $\abs{r} = \sqrt{2}$. Then,
$r = 2\varepsilon_i$, $s = \delta_{k,i-1}$ with $k \leq i-1$.
 \begin{equation*}
 \footnotesize
  \setlength\arraycolsep{1pt}
  \begin{matrix}
  r =  2\varepsilon_i =  & \{0  &  \dots & 0 & 2 & \dots & 2 & 1\}, \\
                         &      &        &   & i &     \\
  \end{matrix}  \qquad
  \begin{matrix}
  s =  \delta_{k,i-1} =  & \{0 &  \dots & 0 & 1  & \dots & 1 & 0 & \dots   & 0 \}. \\
                         &     &        &   & k  &       &   & i    \\
  \end{matrix}
\end{equation*}
Because $r_1$ is a simple root with length $1$, then $r_1 = \delta_{j,j}$.
If $j < i$, then $(r - r_1)_j < 0$, i.e., $r - r_1$ is not a root.
Otherwise, $j \geq i$  and  $(s - r_1)_j < 0$. Here, $s - r_1$ is not a root.
Thus, $r - r_1$ and $s-r_1$ cannot be roots at the same time, i.e., $q$ is a mono-quartet.
For $\abs{s} = \sqrt{2}$ reasoning is the same.
\qed
\medskip

\begin{proposition}[Case $C_n$, $\abs{r_1} = \abs{s_1} = \abs{r} = \abs{s} = 1$]
 \label{eq_necessity}
  Let  $\{r_1, r, s, s_1\}$ be a quartet in $C_n$ such that $\abs{r_1} = \abs{s_1} = \abs{r} = \abs{s} = 1$.
  If both $s - r_1$ and $r - r_1$ are roots, then:

  {\rm(i)} The squares of the lengths of $s - r_1$ and $r - r_1$ are both equal to $1$.

  {\rm(ii)} The inner products of $(s, r_1)$ and $(r, r_1)$ are both equal to $1/2$.

  {\rm(iii)} Roots $r_1$ and $s_1$ are orthogonal.
\end{proposition}

\PerfProof (i)
If $\abs{s - r_1}^2 = 2$ and $\abs{r - r_1}^2 = 2$ then $(s, r_1) = 0$ and $(r, r_1) = 0$.
 Therefore,
 \begin{equation*}
      (r_1 + s_1, r_1) = (r + s, r_1) = 0, \; \text{ i.e., } \; (s_1, r_1) = -1, \\
 \end{equation*}
 which cannot be, since then we would get $\abs{r_1 + s_1}^2 =  2 + 2(s_1, r_1) = 0$.

 Assume, $s - r_1$ and $r - r_1$ are roots with the different lengths:
 \begin{equation}
  \label{eq_diff_lens}
     \abs{s - r_1}^2 = 2, \qquad \abs{r - r_1}^2 = 1
 \end{equation}
 The root $s - r_1$ can be only of type $2\varepsilon_i$, see Table \ref{tab_Cn}. Since $r_1$ is a simple root then
  $r_1 = \delta_{i-1, i-1}$,
  $s = 2\varepsilon_i + \delta_{i-1,i-1} = \gamma_{i-1, i}$.
 \begin{equation*}
 \footnotesize
  \setlength\arraycolsep{1pt}
  \begin{matrix}
  s - r_1 =  2\varepsilon_i =  & \{0  &  \dots & 0 & 2 & \dots & 2 & 1\}, \\
                         &      &        &   & i &     \\
  \end{matrix}  \qquad
  \begin{matrix}
   r_1 =  \delta_{i-1,i-1} =  & \{0 &  \dots & 0 & 1 & 0 & \dots   & 0 \}, \\
                              &     &        &   &   & i    \\
  \end{matrix}  \qquad
  \begin{matrix}
     s =  \gamma_{i-1, i} =  & \{0  &  \dots & 1 & 2 & \dots & 2 & 1\}. \\
                               &      &        &   & i &     \\
  \end{matrix}  \\
\end{equation*}
Now, consider possible cases for $r$. The coordinate $(r)_i$ cannot be non-zero, otherwise  $(s + r)_i \geq 3$.
Thus, $r = \delta_{k, i-1}$, where  $k \leq i-1$, or $r = \delta_{k, i-2}$, where  $k \leq i-2$.
In these cases $r+s = \gamma_{k,i-1}$ or $r+s = \gamma_{k,i}$:
\begin{equation*}
 \footnotesize
  \setlength\arraycolsep{1pt}
 \begin{split}
  & \begin{matrix}
  r = \delta_{k,i-1} =  & \{0 & \dots & 0 & 1 & \dots & 1 & 0 &\dots  &0\}, \\
                        &     &       &   &  k &      &    &  i &     \\
  \end{matrix}  \qquad \;\;
  \begin{matrix}
   r+s =  \gamma_{k,i-1} =  & \{0  &  1 & \dots & 1 & 2 & 2 &  \dots & 2 & 1\}, \quad \text{ or }\\
                            &      &  k &       &   &   & i &  \\
  \end{matrix} \\
  & \begin{matrix}
  r = \delta_{k,i-2} =  & \{0 & \dots & 0 &  1   & \dots & 1  0 &  0 & \dots & 0\}, \\
                        &     &       &   &  k   &       &      &  i         \\
  \end{matrix}  \qquad
  \begin{matrix}
   r+s =  \gamma_{k,i} =  & \{0  &  1 & \dots & 1 & 1 & 2 &  \dots & 2 & 1\}, \\
                            &      &  k &       &   &   & i &  \\
  \end{matrix} \\
 \end{split}
\end{equation*}
So, $r_1 + s_1 = r+s = \gamma_{k,i-1}$ or $\gamma_{k,i}$.
However, in both cases there exist another expansions for $r_1 + s_1$.
In the first case, $r_1 + s_1 = \delta_{k,k} + \gamma_{k+1,i-1}$,
in the second,  $r_1 + s_1 = \delta_{k,k} + \gamma_{k+1,i}$.
In both cases, we get $\delta_{k,k} \prec \delta_{i-1,i-1}$.
Thus, $\{r_1,s_1 \}$ is not an extraspecial pair, and eq. \eqref{eq_diff_lens} cannot happen.
The lengths of $s-r_1$ and $r - r_1$ are both equal to $1$.

(ii) Since $\abs{s - r_1}^2 = 1$, then $2 -2(s,r_1) = 1$, i.e., $(s,r_1) = 1/2$.
Similarly, $(r,r_1) = 1/2$.

(iii) By (ii) we obtain $(r_1+s_1, r_1) = (r+s, r_1) = 1$, i.e., $(s_1, r_1) = 1 - \abs{r_1}^2 = 0$.
\qed

\begin{proposition}[Case $C_n$, $\abs{r_1} = \abs{s_1} = \abs{r} = \abs{s} = 1$.]
 \label{prop_Cn_4}
  Let $\{r_1, r, s, s_1 \}$ be a quartet in $C_n$  such that
  $\abs{r_1} = \abs{s_1} = \abs{r} = \abs{s} = 1$.

  {\rm(i)} Both vectors $s - r_1$ and $r - r_1$ are roots if and only if $(r_1, s_1) = 0$.

  {\rm(ii)} Both vectors $s - r_1$ and $r - r_1$ are roots if and only if $\abs{r_1 + s_1}^2 = 2$.

\end{proposition}
 \PerfProof
 (i)  The necessity of the condition $(r_1, s_1) = 0$ is proved in Proposition \ref{eq_necessity},(iii).
 Let us prove the sufficiency, put $(r_1, s_1) = 0$. Then $(r_1 + s_1, r_1) = (r + s, r_1) = 1$  and, therefore,
 \begin{equation}
  \label{eq_r1_1}
   (r, r_1)+ (s, r_1) = 1.
 \end{equation}
 Since $1 \leq \abs{s - r_1}^2 \leq 2$ then  $1 \leq 2 - 2(s,r_1) \leq  2$. In other words,
 \begin{equation}
  \label{eq_r1_2}
    (s,r_1) = 0 \text{ or } 1/2, \; \text{ similarly  } \;  (r,r_1) = 0 \text{ or } 1/2.
 \end{equation}
 By \eqref{eq_r1_1} and \eqref{eq_r1_2}, we have $(s,r_1) = (r,r_1) = 1/2$.
 Then $\abs{r-r_1}^2 = 1$ and $\abs{s-r_1}^2 = 1$.
 According to Kac's criterion, see Lemma 1.10 b) from \cite{K80},
 we obtain that $r-r_1$ and $s-r_1$ both are roots.

  (ii) This follows from (i) and the following relation:
\begin{equation}
   \abs{r_1 + s_1}^2 = \abs{r_1}^2 + \abs{s_1}^2 + 2(r_1, s_1).
\end{equation}
 \qed

Let us put $\varphi = \displaystyle\frac{\abs{r_1 + s_1}^2}{\abs{s_1}^2}$.

\begin{theorem}
  \label{th_Cn_1}
  {\rm(i)} Let $q = \{r_1, r, s, s_1 \}$ be a quartet in $C_n$.
  The quartet $q$ is a mono-quartet if and only if $(r_1, s_1) \neq 0$.
  This condition does not depend on the pair $\{r, s\}$.

  {\rm(ii)}  The number of extraspecial pairs $\{r_1, s_1\}$ with $(r_1, s_1) = 0$  is $n-1$,
  whereas the total number of extraspecial pairs is $n^2-n$.

   {\rm(iii)} The inner product $(r_1, s_1) = 0$ if and only if $\varphi = 2$.

  {\rm(iv)}  Relation  $\varphi = 1/2$ is true if and only if
  the pair $\{r_1, s_1 \}$ is as follows:
\begin{equation*}
  \footnotesize
  \setlength\arraycolsep{1pt}
   \begin{matrix}
  r_1 =   & \{0 & \dots &  & 1 & 0 &\dots  &0\}, \\
          &     &       &  &   &  i &     \\
  \end{matrix}  \qquad
  \begin{matrix}
   s_1 =   & \{0   & \dots & 0 & 2  &  \dots & 2 & 1\}, \\
           &       &       &   & i  &  \\
  \end{matrix} \qquad
  \begin{matrix}
   r_1+s_1 =  & \{0  & \dots & 1 & 2 &  \dots & 2 & 1\}, \\
              &      &       &   & i &  \\
  \end{matrix}
\end{equation*}
  where $2 \leq i \leq n-1$.
\end{theorem}

 \PerfProof
  (i)  Let $\abs{s_1} = \sqrt{2}$ or $\abs{r_1} = \sqrt{2}$ (case of Proposition \ref{prop_Cn_1}, $q$ is a mono-quartet).
  If $(r_1, s_1) = 0$, we have $\abs{r_1 + s_1}^2 \geq 3$, which cannot be. Therefore, $(r_1, s_1) \neq 0$.
  Let $\abs{s_1} = \abs{r_1} = 1$ and  $\abs{s} = \sqrt{2}$ or $\abs{r} = \sqrt{2}$
  (case of Proposition \ref{prop_Cn_2}, $q$ is a mono-quartet).
  Again, if $(r_1, s_1) = 0$, we have $\abs{r_1 + s_1}^2 = 2$ and $\abs{r + s}^2 \geq 3$, a contradiction.
  So in this case $(r_1, s_1) \neq 0$, too. At last, let
  $\abs{s_1} = \abs{r_1} = \abs{s} = \abs{r} = 1$ (case of Proposition \ref{prop_Cn_4}).
  Here, also, if $q$ is a mono-quartet, then $(r_1, s_1) \neq 0$.
\medskip

  (ii) The required number is
  the number of extraspecial pairs $\{r_1, s_1\}$ with $\abs{r_1 + s_1}^2 = 2$, which is
  the number of the positive roots with length $\sqrt{2}$ that can be decomposed
  into the sum $r_1 + s_1$.  According to Table \ref{tab_Cn}, only roots $2\varepsilon_i$ with $i < n$ satisfy this,
  i.e., the required number is $n-1$.

  By Corollary \ref{col_extrasp}, the total number of extraspecial pairs is $n^2 - n$.
\medskip

  (iii) Relation  $(r_1, s_1) = 0$ is equivalent to
\begin{equation*}
  \abs{r_1 + s_1}^2 = \abs{r_1}^2 + \abs{s_1}^2, \;\text{ i.e., to } \;
  \abs{r_1}^2 = \abs{s_1}^2 = 1 \text { and } \abs{r_1 + s_1}^2 = 2.
\end{equation*}
   Thus, $\varphi = 2$.
\medskip

   (iv) Relation $\varphi = 1/2$  is true if and only if $\abs{r_1 + s_1} = 1$ and $\abs{s_1} = \sqrt{2}$.
Then, according to Table \ref{tab_Cn},
 \begin{equation}
 \footnotesize
 \setlength\arraycolsep{1pt}
 \begin{matrix}
 s_1  = 2\varepsilon_i  =
   & \{0   & \dots & 0 & 2  &  \dots & 2 & 1\}. \\
   &       &       &   & i  &  \\
  \end{matrix}
\end{equation}
 Since $r_1$ is a simple root, by Table \ref{tab_Cn} there is only one possibility to $r_1$:
 \begin{equation}
 \footnotesize
 \setlength\arraycolsep{1pt}
    \begin{matrix}
  r_1 =  \delta_{i-1, i-1}  =
      & \{0 & \dots &  & 1 & 0 &\dots  &0\}, \\
      &     &       &  &   &  i &     \\
  \end{matrix}
\end{equation}
\qed

  Similarly to  Theorem \ref{th_Bn_2} for $B_n$, the following theorem gives \ref{th_Bn_1} more precise relation
  for the lengths $s - r_1$ and $r - r_1$.

\begin{theorem}
  \label{th_Cn_2}
  In the root system $C_n$, any quartet is simple.
\end{theorem}

\PerfProof  Since $r_1$ is a simple root,
either $r_1 = \delta_{i,i}$ with $i < n$, or $r_1 = 2\varepsilon_n$, see Table \ref{tab_Cn}.
For $\abs{r_1} =  2\varepsilon_n$ (the only simple root with length $\sqrt{2}$)
or $\abs{s_1} = \sqrt{2}$, by Proposition \ref{prop_Cn_1} there is no quartet.
Thus, we can assume that $\abs{r_1} = \abs{s_1} = 1$.

Consider the case when $s-r_1$ is a root. The case when $r-r_1$ is a root is treated similarly.
By Proposition \ref{prop_Cn_2}, if $s-r_1$ is a root then $\abs{s}$ cannot be of length $\sqrt{2}$.
Therefore, $\abs{s} = 1$. We will show that $\abs{s-r_1} = 1$, i.e., $\abs{s-r_1} = \abs{s}$.

Since $r_1$ is a simple root of length $1$, $r_1 = \delta_{i,i}$.
We assume that $s - r_1$ is a root of length $\sqrt{2}$, so $s - r_1 = 2\varepsilon_{i+1}$, and
$s = \delta_{i,i} + 2\varepsilon_{i+1}$:
\begin{equation}
 \label{eq_lens_eq_1}
 \footnotesize
  \setlength\arraycolsep{1pt}
  \begin{matrix}
  r_1 = \delta_{i,i} =  & \{0 & \dots & 0 & 1 & 0 & \dots &  & 0 \}, \\
                        &     &       &   & i &      &    &       \\
  \end{matrix}  \quad
  \begin{matrix}
  s - r_1 =  2\varepsilon_{i+1} =  & \{0  &  \dots  & 0  & 2 & 2 &  \dots & 2 & 1\}, \\
                               &      &         & i  &       \\
  \end{matrix}
  \begin{matrix} \quad
  s = \gamma_{i,i+1} =   & \{0  &  \dots  & 1  & 2 & 2 &  \dots & 2 & 1\}. \\
                       &      &         & i  &       \\
  \end{matrix}
\end{equation}
For $s$ in \eqref{eq_lens_eq_1}, there are only two options for $r$ such that $r+s$ is the root:
\begin{equation}
 \label{eq_lens_eq_2}
 \footnotesize
  \setlength\arraycolsep{1pt}
 \begin{split}
 & \begin{matrix}
   \qquad  r =  \delta_{k,i-1} =    & \{0  &  \dots & 0 & 1 & \dots & 1  &  0 & \dots 0 \}, \\
                                    &      &        &   & k &       &    &  i  &         \\
  \end{matrix} \qquad
  \begin{matrix}
   r + s =  & \gamma_{k,i+1} =     & \{0  &  \dots 0 & 1 & \dots    &  1 &  2 \dots & 2 & 1\}, \text{ or} \\
            &                          &      &        & k &          &  i \\
  \end{matrix} \\
 & \begin{matrix}
   \qquad  r =  & \delta_{k,i} \quad =    & \{0  \dots 0 & 1 & \dots  & 1  &  0 & \dots & 0 \}, \\
                &                  &              & k &        & i  \\
  \end{matrix} \qquad
  \begin{matrix}
   r + s =  \gamma_{k,i} =     & \{0  &  \dots & 0 & 1 & \dots & 1 &  2 & \dots 2 & 1\}.\\
                               &      &        &   & k &       &   &  i & \\
  \end{matrix} \\
 \end{split}
\end{equation}
 In both cases \eqref{eq_lens_eq_2} $k < i$, and for both there is another expansion for $r_1 + s_1 = r + s$
 with $r = \delta_{k,k}$. Since $r_1 = \delta_{i,i}$, then $r \prec r_1$, which cannot be.  \qed
~\\

It is easy to check the following relations for the parameter $\varphi$:
\begin{equation}
  \label{eq_cases_fi}
    \varphi =
 \begin{cases}
  \footnotesize
  \setlength\arraycolsep{1pt}
   \begin{split}
    2  \; \Longleftrightarrow \; & \abs{r_1 + s_1} = \sqrt{2} \text{ and }\abs{s_1} = 1, \\
        &  \begin{matrix}
               r_1 + s_1 = 2\varepsilon_i = & \{0  &  \dots  & 0  & 2 &  \dots & 2 & 1\},
                    \text{ where } 1 \leq i \leq n-1,  \\
                                      &      &         &    & i  &       \\
            \end{matrix}
   \end{split}  \\
   \footnotesize
   \setlength\arraycolsep{1pt}
   \begin{split}
   1/2 \; \Longleftrightarrow \; & \abs{r_1 + s_1} = 1 \text{ and }\abs{s_1} = \sqrt{2}, \\
         & \begin{matrix}
               s_1 = 2\varepsilon_i = & \{0  &  \dots  & 0  & 2 &  \dots & 2 & 1\},
                    \text{ where } 1 \leq i \leq n,  \\
                                      &      &         &    & i  &       \\
            \end{matrix}
   \end{split}  \\
   \footnotesize
   \begin{split}
    1 \; \Longleftrightarrow \; & \abs{r_1 + s_1} = 1 \text{ and }\abs{s_1} = 1, \\
         & \text{ The case } \abs{r_1 + s_1} = \sqrt{2} \text{ and }\abs{s_1} = \sqrt{2} \text{ cannot happen.}
   \end{split}  \\
 \end{cases}
\end{equation} \qed

\subsection{Extraspecial pairs and quartets in $F_4$}
  \label{sec_F4}

Although the quartets in $F_4$ are not structured as well as in $B_n$ and $C_n$,
we can find all simple quartets.
In the case $F_4$ there are $48$ quartets, of which $38$ have
the same good properties as in $B_n$ and $C_n$, however, there  are $10$ quartets being an obstacle
to using formulas \eqref{eq_Cn_final_0} or \eqref{eq_Bn_final_0}.
 The positive roots of $F_4$, ordered in the regular ordering,
are shown in eq.~\eqref{eq_F4_roots}.

\begin{equation}
 \label{eq_F4_roots}
  \begin{matrix}
0 )  & [1, 0, 0, 0]  &   2 \\
1 )  & [0, 1, 0, 0]  &   2 \\
2 )  & [0, 0, 1, 0]  &   1 \\
3 )  & [0, 0, 0, 1]  &   1 \\
4 )  & [1, 1, 0, 0]  &   2 \\
5 )  & [0, 1, 1, 0]  &   1 \\
  \end{matrix}
  \qquad
  \begin{matrix}
6 )  & [0, 0, 1, 1]  &   1 \\
7 )  & [1, 1, 1, 0]  &   1 \\
8 )  & [0, 1, 2, 0]  &   2 \\
9 )  & [0, 1, 1, 1]  &   1 \\
10 )  & [1, 1, 2, 0]  &   2 \\
11 )  & [1, 1, 1, 1]  &   1 \\
  \end{matrix}
  \qquad
  \begin{matrix}
12 )  & [0, 1, 2, 1]  &   1 \\
13 )  & [1, 2, 2, 0]  &   2 \\
14 )  & [1, 1, 2, 1]  &   1 \\
15 )  & [0, 1, 2, 2]  &   2 \\
16 )  & [1, 2, 2, 1]  &   1 \\
17 )  & [1, 1, 2, 2]  &   2 \\
  \end{matrix}
  \qquad
  \begin{matrix}
18 )  & [1, 2, 3, 1]  &   1 \\
19 )  & [1, 2, 2, 2]  &   2 \\
20 )  & [1, 2, 3, 2]  &   1 \\
21 )  & [1, 2, 4, 2]  &   2 \\
22 )  & [1, 3, 4, 2]  &   2 \\
23 )  & [2, 3, 4, 2]  &   2 \\
  \end{matrix}
\end{equation}

\begin{proposition}[Case $F_4$, $\abs{r_1} = \sqrt{2}$]
  \label{prop_F4_1}
 {\rm(i)} Let $\alpha_i$ be the root with index $i$ of the ordered list in eq.~\eqref{eq_F4_roots},
 and $\theta$ be an arbitrary positive root in $F_4$.

 {\rm(i)} If $\alpha$ is a simple root of length $\sqrt{2}$,
 the inner product $(\theta, \alpha)$ is an integer number.

 {\rm(ii)} For an extraspecial root $\{r_1, s_1 \}$ with $\abs{r_1} = \sqrt{2}$, we have $(r_1, s_1) = -1$.

{\rm(iii)}   Let $q = \{r_1, r, s, s_1 \}$ be a quartet with $\abs{r_1} = \sqrt{2}$.
  Then $q$ is a mono-quartet.
\end{proposition}

\PerfProof (i) The inner products of roots in $F_4$ are $0, \pm{1/2}, \pm{1}$.
However, the inner products for  $(\theta, \alpha)$, where $\alpha = \alpha_i$, $i = 1,2$ can be only
$0, \pm{1}$. To see this, let us get the matrix $B$ of the bilinear form for $F_4$:
\begin{equation}
 \label{eq_F4_bil_form}
  B = \begin{bmatrix}
        2 & -1 &  0 & 0 \\
       -1 &  2 & -1 & 0 \\
        0 & -1 &  1 & -1/2 \\
        0 &  0 & -1/2 & 1 \\
      \end{bmatrix}
\end{equation}
The inner products  $(\theta, \alpha_1)$ and $(\theta, \alpha_2)$ are as follows:
\begin{equation*}
  (\theta, \alpha_1) = \theta^t B \alpha_1 =  2\theta_1 - \theta_2,
  \qquad (\theta, \alpha_2) = \theta^t B \alpha_2 = -\theta_1 + 2\theta_2 - \theta_3.
\end{equation*}

Obviously, that equations  $2\theta_1 - \theta_2 = \pm{1/2}$ and
$- \theta_1 + 2\theta_2 - \theta_3 = \pm{1/2}$  have no solutions
for any integer vector $\theta$.

 The proof of (ii) and (iii) coincides with the proof of Proposition \ref{prop_Bn_1}(i),(ii)
for the case $B_n$. \qed
\bigskip
\bigskip

\begin{proposition}[Case $F_4$, $\abs{r_1} = 1$]
   Let $q = \{r_1, r, s, s_1 \}$ be a simple quartet with
\begin{equation*}
   \abs{r_1} = 1 \text{ and } \abs{r} = 1 \quad (\text{resp. } \abs{s} = 1).
\end{equation*}
    Then, $r_1 = \alpha_2$ or $r_1 = \alpha_3$ and the following properties hold:

   {\rm(i)}  $(r,r_1) = 1/2$  (resp. $(s,r_1) = 1/2$).

   {\rm(ii)} For $r_1 = \alpha_2$,  coordinates of $r$ (resp. $s$) are as follows:
\begin{equation}
  \label{eq_cond_F4_1}
       (r)_4 = 1,\, (r)_3 = (r)_2 + 1  \quad (\text{resp. } (s)_4 = 1,\, (s)_3 = (s)_2 + 1).
\end{equation}

   {\rm(iii)} For $r_1 = \alpha_3$,  coordinates of $r$ (resp. $s$) are as follows:
\begin{equation}
  \label{eq_cond_F4_2}
       (r)_3 + 1 = 2(r)_4 \quad (\text{resp. } (s)_3 + 1 = 2(s)_4).
\end{equation}
\end{proposition}

\PerfProof (i) Consider case $\abs{r} = 1$. Here,
\begin{equation*}
   \abs{r}^2 = \abs{r - r_1}^2 \; \Longleftrightarrow \; 2 - 2(r, r_1) = 1 \; \Longleftrightarrow \; (r, r_1) = 1/2.
\end{equation*}

(ii) Since $r_1$ is a simple root with $\abs{r_1} = 1$, then $r_1 = \alpha_2$ or $r_1 = \alpha_3$.
For $r_1 = \alpha_2$, by \eqref{eq_F4_bil_form}, we have
\begin{equation*}
 \begin{split}
  & (r,r_1) = (r, \alpha_2) = r^t B \alpha_2 = -(r)_2 + (r)_3 - (r_4)/2  = 1/2, \; \text{ i.e., }  \\
  & (r)_3 - (r)_2 = ((r)_4 + 1)/2.
 \end{split}
\end{equation*}
Since $(r)_4$ takes values $0$, $1$ or $2$, see eq.~\eqref{eq_F4_roots}, we have  $(r)_4 = 1$ and  $(r)_3 = (r)_2 + 1$.

(iii) For $r_1 = \alpha_3$, by \eqref{eq_F4_bil_form},
\begin{equation*}
  (r,r_1) = (r, \alpha_3) = r^t B \alpha_3 = -(r)_3/2 + (r)_4 = 1/2.
\end{equation*}
Therefore, $(r)_3 + 1 = 2(r)_4$. \qed

\begin{corollary}
 \label{col_F4_1}
  {\rm(i)} Only roots $\alpha_6$, $\alpha_{12}$, $\alpha_{14}$, $\alpha_{18}$ satisfy eq. \eqref{eq_cond_F4_1}.
     Roots $\alpha_6$ and $\alpha_9$ satisfy eq. \eqref{eq_cond_F4_2}.

  {\rm(ii)} The following  quartets with $\abs{r_1} = 1$ are simple:
\begin{equation}
 \label{eq_F4_8_quart}
\begin{array}{llllllll}
  & [r_1, r, s, s_1] = [ 2 &  \fbox{6} &  13 &  16 ] & \quad r_1 = \alpha_2  & \quad \abs{r-r_1} = \abs{r}  & \\
  & [r_1, r, s, s_1] = [ 2 &  3 &  \fbox{18} &  19 ] & \quad r_1 = \alpha_2  & \quad \abs{s-r_1} = \abs{s}  &  \\
  & [r_1, r, s, s_1] = [ 2 &  \fbox{6} &  16 &  19 ] & \quad r_1 = \alpha_2  & \quad \abs{r-r_1} = \abs{r}  & \\
  & [r_1, r, s, s_1] = [ 2 &  9 &  \fbox{14} &  19 ] & \quad r_1 = \alpha_2  & \quad \abs{s-r_1} = \abs{s}  &   \\
  & [r_1, r, s, s_1] = [ 2 &  11 &  \fbox{12} &  19 ] & \quad r_1 = \alpha_2 & \quad \abs{s-r_1} = \abs{s}  &   \\
  & [r_1, r, s, s_1] = [ 2 &  \fbox{6} &  \fbox{18} &  20 ] & \quad r_1 = \alpha_2  & \quad \abs{s-r_1} = \abs{s},  & \abs{r-r_1} = \abs{r} \\
  & [r_1, r, s, s_1] = [ 2 &  \fbox{12} &  \fbox{14} &  20 ] & \quad r_1 = \alpha_2 & \quad \abs{s-r_1} = \abs{s},  & \abs{r-r_1} = \abs{r}  \\
  & [r_1, r, s, s_1] = [ 3 &  \fbox{6} &  \fbox{9} &  12 ] & \quad r_1 = \alpha_3   & \quad \abs{s-r_1} = \abs{s},  & \abs{r-r_1} = \abs{r}  \\
  \end{array}
\end{equation}
  The roots enclosed in the box are the roots from (i).
\end{corollary}
\qed

\begin{theorem}
 \label{th_F4_1}
In the root system $F_4$,

{\rm(i)}  There are $30$ simple quartets with $r_1 = \alpha_0$ or $r_1 = \alpha_1$.
see Table~\ref{table_F4_quartets}.
 If $q = \{r_1, r, s, s_1 \}$ such a quartet, then
 roots $r + s$ and $s_1$ have the same length.

{\rm(ii)} There are $8$ simple quartets with $r_1 = \alpha_2$ or $r_1 = \alpha_3$ such that
$\abs{r} = 1$ or $\abs{s} = 1$. They listed in eq.~\eqref{eq_F4_8_quart}.

{\rm(iii)}  The remaining $10$ quartets are not simple. These quartets are as follows:
$31-33$, $35-37$, $39$, $41$, $45$, $46$ see Table~\ref{table_F4_quartets}.
\end{theorem}

\PerfProof
(i) If $s - r_1$ is a root, then $\abs{s - r_1}^2$ = $\abs{s}^2 + 2 - 2(s, r_1) \leq 2$, i.e.,
$(s, r_1) > 0$.   Since $\abs{r_1} = \sqrt{2}$, by Proposition \ref{prop_F4_1}(i),
the inner product $(s, r_1)$ cannot be $1/2$.
Then, $(s, r_1) = 1$, and $\abs{s - r_1}^2 = \abs{s}^2$.
Similarly, if $r - r_1$ is a root, then $\abs{r - r_1}^2 = \abs{r}^2$.
Since $r + s = r_1 + s_1$,
by Proposition \ref{prop_F4_1}(ii),
the length of $r + s$ is equal to $2 + 2(r_1, s_1) + \abs{s_1}^2 = \abs{s_1}^2$.

(ii) Follows from Corollary \ref{col_F4_1}.

 (iii) This  verified simply by calculating the lengths for the corresponding roots.
\qed

\begin{appendix}
\clearpage
\section{\bf The quartets for $B_6$, $C_6$ and $F_4$}
  \label{sec_app}

\begin{table}[!ht]
$\begin{matrix}
0 )  & [1, 0, 0, 0, 0, 0]  &   2 \\
1 )  & [0, 1, 0, 0, 0, 0]  &   2 \\
2 )  & [0, 0, 1, 0, 0, 0]  &   2 \\
3 )  & [0, 0, 0, 1, 0, 0]  &   2 \\
4 )  & [0, 0, 0, 0, 1, 0]  &   2 \\
5 )  & [0, 0, 0, 0, 0, 1]  &   1 \\
6 )  & [1, 1, 0, 0, 0, 0]  &   2 \\
7 )  & [0, 1, 1, 0, 0, 0]  &   2 \\
8 )  & [0, 0, 1, 1, 0, 0]  &   2 \\
9 )  & [0, 0, 0, 1, 1, 0]  &   2 \\
10 )  & [0, 0, 0, 0, 1, 1]  &   1 \\
11 )  & [1, 1, 1, 0, 0, 0]  &   2 \\
12 )  & [0, 1, 1, 1, 0, 0]  &   2 \\
13 )  & [0, 0, 1, 1, 1, 0]  &   2 \\
14 )  & [0, 0, 0, 1, 1, 1]  &   1 \\
15 )  & [0, 0, 0, 0, 1, 2]  &   2 \\
16 )  & [1, 1, 1, 1, 0, 0]  &   2 \\
17 )  & [0, 1, 1, 1, 1, 0]  &   2 \\
18 )  & [0, 0, 1, 1, 1, 1]  &   1 \\
19 )  & [0, 0, 0, 1, 1, 2]  &   2 \\
20 )  & [1, 1, 1, 1, 1, 0]  &   2 \\
21 )  & [0, 1, 1, 1, 1, 1]  &   1 \\
22 )  & [0, 0, 1, 1, 1, 2]  &   2 \\
23 )  & [0, 0, 0, 1, 2, 2]  &   2 \\
24 )  & [1, 1, 1, 1, 1, 1]  &   1 \\
25 )  & [0, 1, 1, 1, 1, 2]  &   2 \\
26 )  & [0, 0, 1, 1, 2, 2]  &   2 \\
27 )  & [1, 1, 1, 1, 1, 2]  &   2 \\
28 )  & [0, 1, 1, 1, 2, 2]  &   2 \\
29 )  & [0, 0, 1, 2, 2, 2]  &   2 \\
30 )  & [1, 1, 1, 1, 2, 2]  &   2 \\
31 )  & [0, 1, 1, 2, 2, 2]  &   2 \\
32 )  & [1, 1, 1, 2, 2, 2]  &   2 \\
33 )  & [0, 1, 2, 2, 2, 2]  &   2 \\
34 )  & [1, 1, 2, 2, 2, 2]  &   2 \\
35 )  & [1, 2, 2, 2, 2, 2]  &   2 \\
\end{matrix}$
 \qquad \qquad
 \qquad \qquad
$\begin{matrix}
0 )  & [1, 0, 0, 0, 0, 0]  &   1 \\
1 )  & [0, 1, 0, 0, 0, 0]  &   1 \\
2 )  & [0, 0, 1, 0, 0, 0]  &   1 \\
3 )  & [0, 0, 0, 1, 0, 0]  &   1 \\
4 )  & [0, 0, 0, 0, 1, 0]  &   1 \\
5 )  & [0, 0, 0, 0, 0, 1]  &   2 \\
6 )  & [1, 1, 0, 0, 0, 0]  &   1 \\
7 )  & [0, 1, 1, 0, 0, 0]  &   1 \\
8 )  & [0, 0, 1, 1, 0, 0]  &   1 \\
9 )  & [0, 0, 0, 1, 1, 0]  &   1 \\
10 )  & [0, 0, 0, 0, 1, 1]  &   1 \\
11 )  & [1, 1, 1, 0, 0, 0]  &   1 \\
12 )  & [0, 1, 1, 1, 0, 0]  &   1 \\
13 )  & [0, 0, 1, 1, 1, 0]  &   1 \\
14 )  & [0, 0, 0, 1, 1, 1]  &   1 \\
15 )  & [0, 0, 0, 0, 2, 1]  &   2 \\
16 )  & [1, 1, 1, 1, 0, 0]  &   1 \\
17 )  & [0, 1, 1, 1, 1, 0]  &   1 \\
18 )  & [0, 0, 1, 1, 1, 1]  &   1 \\
19 )  & [0, 0, 0, 1, 2, 1]  &   1 \\
20 )  & [1, 1, 1, 1, 1, 0]  &   1 \\
21 )  & [0, 1, 1, 1, 1, 1]  &   1 \\
22 )  & [0, 0, 1, 1, 2, 1]  &   1 \\
23 )  & [0, 0, 0, 2, 2, 1]  &   2 \\
24 )  & [1, 1, 1, 1, 1, 1]  &   1 \\
25 )  & [0, 1, 1, 1, 2, 1]  &   1 \\
26 )  & [0, 0, 1, 2, 2, 1]  &   1 \\
27 )  & [1, 1, 1, 1, 2, 1]  &   1 \\
28 )  & [0, 1, 1, 2, 2, 1]  &   1 \\
29 )  & [0, 0, 2, 2, 2, 1]  &   2 \\
30 )  & [1, 1, 1, 2, 2, 1]  &   1 \\
31 )  & [0, 1, 2, 2, 2, 1]  &   1 \\
32 )  & [1, 1, 2, 2, 2, 1]  &   1 \\
33 )  & [0, 2, 2, 2, 2, 1]  &   2 \\
34 )  & [1, 2, 2, 2, 2, 1]  &   1 \\
35 )  & [2, 2, 2, 2, 2, 1]  &   2 \\
\end{matrix}$
\medskip
\caption{ The regular ordering of the positive roots in $B_6$ (left) and $C_6$ (right).
The squares of the root lengths are also indicated.}
\label{tab_B6C6}
\end{table}

\begin{table}[!ht]
$ \qquad
  \tiny
  \begin{array}{llllllll}
1 ) &  [r_1, r, s, s_1] = [ 0 &  2 &  6 &  7 ] & \quad r + s = \alpha_{11}  & \quad s-r_1 = \alpha_{ 1 }  &   \\
2 ) &  [r_1, r, s, s_1] = [ 0 &  3 &  11 &  12 ] & \quad r + s = \alpha_{16}  & \quad s-r_1 = \alpha_{ 7 }  &   \\
3 ) &  [r_1, r, s, s_1] = [ 0 &  6 &  8 &  12 ] & \quad r + s = \alpha_{16}  &   & \quad r-r_1 = \alpha_{ 1 }  \\
4 ) &  [r_1, r, s, s_1] = [ 0 &  4 &  16 &  17 ] & \quad r + s = \alpha_{20}  & \quad s-r_1 = \alpha_{ 12 }  &   \\
5 ) &  [r_1, r, s, s_1] = [ 0 &  6 &  13 &  17 ] & \quad r + s = \alpha_{20}  &   & \quad r-r_1 = \alpha_{ 1 }  \\
6 ) &  [r_1, r, s, s_1] = [ 0 &  9 &  11 &  17 ] & \quad r + s = \alpha_{20}  & \quad s-r_1 = \alpha_{ 7 }  &   \\
7 ) &  [r_1, r, s, s_1] = [ 0 &  5 &  20 &  21 ] & \quad r + s = \alpha_{24}  & \quad s-r_1 = \alpha_{ 17 }  &   \\
8 ) &  [r_1, r, s, s_1] = [ 0 &  6 &  18 &  21 ] & \quad r + s = \alpha_{24}  &   & \quad r-r_1 = \alpha_{ 1 }  \\
9 ) &  [r_1, r, s, s_1] = [ 0 &  10 &  16 &  21 ] & \quad r + s = \alpha_{24}  & \quad s-r_1 = \alpha_{ 12 }  &   \\
10 ) &  [r_1, r, s, s_1] = [ 0 &  11 &  14 &  21 ] & \quad r + s = \alpha_{24}  &   & \quad r-r_1 = \alpha_{ 7 }  \\
11 ) &  [r_1, r, s, s_1] = [ 0 &  5 &  24 &  25 ] & \quad r + s = \alpha_{27}  & \quad s-r_1 = \alpha_{ 21 }  &   \\
12 ) &  [r_1, r, s, s_1] = [ 0 &  6 &  22 &  25 ] & \quad r + s = \alpha_{27}  &   & \quad r-r_1 = \alpha_{ 1 }  \\
13 ) &  [r_1, r, s, s_1] = [ 0 &  11 &  19 &  25 ] & \quad r + s = \alpha_{27}  &   & \quad r-r_1 = \alpha_{ 7 }  \\
14 ) &  [r_1, r, s, s_1] = [ 0 &  15 &  16 &  25 ] & \quad r + s = \alpha_{27}  & \quad s-r_1 = \alpha_{ 12 }  &   \\
15 ) &  [r_1, r, s, s_1] = [ 0 &  4 &  27 &  28 ] & \quad r + s = \alpha_{30}  & \quad s-r_1 = \alpha_{ 25 }  &   \\
16 ) &  [r_1, r, s, s_1] = [ 0 &  6 &  26 &  28 ] & \quad r + s = \alpha_{30}  &   & \quad r-r_1 = \alpha_{ 1 }  \\
17 ) &  [r_1, r, s, s_1] = [ 0 &  10 &  24 &  28 ] & \quad r + s = \alpha_{30}  & \quad s-r_1 = \alpha_{ 21 }  &   \\
18 ) &  [r_1, r, s, s_1] = [ 0 &  11 &  23 &  28 ] & \quad r + s = \alpha_{30}  &   & \quad r-r_1 = \alpha_{ 7 }  \\
19 ) &  [r_1, r, s, s_1] = [ 0 &  15 &  20 &  28 ] & \quad r + s = \alpha_{30}  & \quad s-r_1 = \alpha_{ 17 }  &   \\
20 ) &  [r_1, r, s, s_1] = [ 0 &  3 &  30 &  31 ] & \quad r + s = \alpha_{32}  & \quad s-r_1 = \alpha_{ 28 }  &   \\
21 ) &  [r_1, r, s, s_1] = [ 0 &  6 &  29 &  31 ] & \quad r + s = \alpha_{32}  &   & \quad r-r_1 = \alpha_{ 1 }  \\
22 ) &  [r_1, r, s, s_1] = [ 0 &  9 &  27 &  31 ] & \quad r + s = \alpha_{32}  & \quad s-r_1 = \alpha_{ 25 }  &   \\
23 ) &  [r_1, r, s, s_1] = [ 0 &  14 &  24 &  31 ] & \quad r + s = \alpha_{32}  & \quad s-r_1 = \alpha_{ 21 }  &   \\
24 ) &  [r_1, r, s, s_1] = [ 0 &  16 &  23 &  31 ] & \quad r + s = \alpha_{32}  &   & \quad r-r_1 = \alpha_{ 12 }  \\
25 ) &  [r_1, r, s, s_1] = [ 0 &  19 &  20 &  31 ] & \quad r + s = \alpha_{32}  & \quad s-r_1 = \alpha_{ 17 }  &   \\
26 ) &  [r_1, r, s, s_1] = [ 0 &  2 &  32 &  33 ] & \quad r + s = \alpha_{34}  & \quad s-r_1 = \alpha_{ 31 }  &   \\
27 ) &  [r_1, r, s, s_1] = [ 0 &  8 &  30 &  33 ] & \quad r + s = \alpha_{34}  & \quad s-r_1 = \alpha_{ 28 }  &   \\
28 ) &  [r_1, r, s, s_1] = [ 0 &  11 &  29 &  33 ] & \quad r + s = \alpha_{34}  &   & \quad r-r_1 = \alpha_{ 7 }  \\
29 ) &  [r_1, r, s, s_1] = [ 0 &  13 &  27 &  33 ] & \quad r + s = \alpha_{34}  & \quad s-r_1 = \alpha_{ 25 }  &   \\
30 ) &  [r_1, r, s, s_1] = [ 0 &  16 &  26 &  33 ] & \quad r + s = \alpha_{34}  &   & \quad r-r_1 = \alpha_{ 12 }  \\
31 ) &  [r_1, r, s, s_1] = [ 0 &  18 &  24 &  33 ] & \quad r + s = \alpha_{34}  & \quad s-r_1 = \alpha_{ 21 }  &   \\
32 ) &  [r_1, r, s, s_1] = [ 0 &  20 &  22 &  33 ] & \quad r + s = \alpha_{34}  &   & \quad r-r_1 = \alpha_{ 17 }  \\
33 ) &  [r_1, r, s, s_1] = [ 1 &  3 &  7 &  8 ] & \quad r + s = \alpha_{12}  & \quad s-r_1 = \alpha_{ 2 }  &   \\
34 ) &  [r_1, r, s, s_1] = [ 1 &  4 &  12 &  13 ] & \quad r + s = \alpha_{17}  & \quad s-r_1 = \alpha_{ 8 }  &   \\
35 ) &  [r_1, r, s, s_1] = [ 1 &  7 &  9 &  13 ] & \quad r + s = \alpha_{17}  &   & \quad r-r_1 = \alpha_{ 2 }  \\
36 ) &  [r_1, r, s, s_1] = [ 1 &  5 &  17 &  18 ] & \quad r + s = \alpha_{21}  & \quad s-r_1 = \alpha_{ 13 }  &   \\
37 ) &  [r_1, r, s, s_1] = [ 1 &  7 &  14 &  18 ] & \quad r + s = \alpha_{21}  &   & \quad r-r_1 = \alpha_{ 2 }  \\
38 ) &  [r_1, r, s, s_1] = [ 1 &  10 &  12 &  18 ] & \quad r + s = \alpha_{21}  & \quad s-r_1 = \alpha_{ 8 }  &   \\
39 ) &  [r_1, r, s, s_1] = [ 1 &  5 &  21 &  22 ] & \quad r + s = \alpha_{25}  & \quad s-r_1 = \alpha_{ 18 }  &   \\
40 ) &  [r_1, r, s, s_1] = [ 1 &  7 &  19 &  22 ] & \quad r + s = \alpha_{25}  &   & \quad r-r_1 = \alpha_{ 2 }  \\
41 ) &  [r_1, r, s, s_1] = [ 1 &  12 &  15 &  22 ] & \quad r + s = \alpha_{25}  &   & \quad r-r_1 = \alpha_{ 8 }  \\
42 ) &  [r_1, r, s, s_1] = [ 1 &  4 &  25 &  26 ] & \quad r + s = \alpha_{28}  & \quad s-r_1 = \alpha_{ 22 }  &   \\
43 ) &  [r_1, r, s, s_1] = [ 1 &  7 &  23 &  26 ] & \quad r + s = \alpha_{28}  &   & \quad r-r_1 = \alpha_{ 2 }  \\
44 ) &  [r_1, r, s, s_1] = [ 1 &  10 &  21 &  26 ] & \quad r + s = \alpha_{28}  & \quad s-r_1 = \alpha_{ 18 }  &   \\
45 ) &  [r_1, r, s, s_1] = [ 1 &  15 &  17 &  26 ] & \quad r + s = \alpha_{28}  & \quad s-r_1 = \alpha_{ 13 }  &   \\
46 ) &  [r_1, r, s, s_1] = [ 1 &  3 &  28 &  29 ] & \quad r + s = \alpha_{31}  & \quad s-r_1 = \alpha_{ 26 }  &   \\
47 ) &  [r_1, r, s, s_1] = [ 1 &  9 &  25 &  29 ] & \quad r + s = \alpha_{31}  & \quad s-r_1 = \alpha_{ 22 }  &   \\
48 ) &  [r_1, r, s, s_1] = [ 1 &  12 &  23 &  29 ] & \quad r + s = \alpha_{31}  &   & \quad r-r_1 = \alpha_{ 8 }  \\
49 ) &  [r_1, r, s, s_1] = [ 1 &  14 &  21 &  29 ] & \quad r + s = \alpha_{31}  & \quad s-r_1 = \alpha_{ 18 }  &   \\
50 ) &  [r_1, r, s, s_1] = [ 1 &  17 &  19 &  29 ] & \quad r + s = \alpha_{31}  &   & \quad r-r_1 = \alpha_{ 13 }  \\
51 ) &  [r_1, r, s, s_1] = [ 1 &  6 &  33 &  34 ] & \quad r + s = \alpha_{35}  &   & \quad r-r_1 = \alpha_{ 0 }  \\
52 ) &  [r_1, r, s, s_1] = [ 1 &  7 &  32 &  34 ] & \quad r + s = \alpha_{35}  &   & \quad r-r_1 = \alpha_{ 2 }  \\
53 ) &  [r_1, r, s, s_1] = [ 1 &  11 &  31 &  34 ] & \quad r + s = \alpha_{35}  & \quad s-r_1 = \alpha_{ 29 }  &   \\
54 ) &  [r_1, r, s, s_1] = [ 1 &  12 &  30 &  34 ] & \quad r + s = \alpha_{35}  &   & \quad r-r_1 = \alpha_{ 8 }  \\
55 ) &  [r_1, r, s, s_1] = [ 1 &  16 &  28 &  34 ] & \quad r + s = \alpha_{35}  & \quad s-r_1 = \alpha_{ 26 }  &   \\
56 ) &  [r_1, r, s, s_1] = [ 1 &  17 &  27 &  34 ] & \quad r + s = \alpha_{35}  &   & \quad r-r_1 = \alpha_{ 13 }  \\
57 ) &  [r_1, r, s, s_1] = [ 1 &  20 &  25 &  34 ] & \quad r + s = \alpha_{35}  & \quad s-r_1 = \alpha_{ 22 }  &   \\
58 ) &  [r_1, r, s, s_1] = [ 1 &  21 &  24 &  34 ] & \quad r + s = \alpha_{35}  &   & \quad r-r_1 = \alpha_{ 18 }  \\
59 ) &  [r_1, r, s, s_1] = [ 2 &  4 &  8 &  9 ] & \quad r + s = \alpha_{13}  & \quad s-r_1 = \alpha_{ 3 }  &   \\
60 ) &  [r_1, r, s, s_1] = [ 2 &  5 &  13 &  14 ] & \quad r + s = \alpha_{18}  & \quad s-r_1 = \alpha_{ 9 }  &   \\
61 ) &  [r_1, r, s, s_1] = [ 2 &  8 &  10 &  14 ] & \quad r + s = \alpha_{18}  &   & \quad r-r_1 = \alpha_{ 3 }  \\
62 ) &  [r_1, r, s, s_1] = [ 2 &  5 &  18 &  19 ] & \quad r + s = \alpha_{22}  & \quad s-r_1 = \alpha_{ 14 }  &   \\
63 ) &  [r_1, r, s, s_1] = [ 2 &  8 &  15 &  19 ] & \quad r + s = \alpha_{22}  &   & \quad r-r_1 = \alpha_{ 3 }  \\
64 ) &  [r_1, r, s, s_1] = [ 2 &  4 &  22 &  23 ] & \quad r + s = \alpha_{26}  & \quad s-r_1 = \alpha_{ 19 }  &   \\
65 ) &  [r_1, r, s, s_1] = [ 2 &  10 &  18 &  23 ] & \quad r + s = \alpha_{26}  & \quad s-r_1 = \alpha_{ 14 }  &   \\
66 ) &  [r_1, r, s, s_1] = [ 2 &  13 &  15 &  23 ] & \quad r + s = \alpha_{26}  &   & \quad r-r_1 = \alpha_{ 9 }  \\
67 ) &  [r_1, r, s, s_1] = [ 2 &  7 &  29 &  31 ] & \quad r + s = \alpha_{33}  &   & \quad r-r_1 = \alpha_{ 1 }  \\
68 ) &  [r_1, r, s, s_1] = [ 2 &  8 &  28 &  31 ] & \quad r + s = \alpha_{33}  &   & \quad r-r_1 = \alpha_{ 3 }  \\
69 ) &  [r_1, r, s, s_1] = [ 2 &  12 &  26 &  31 ] & \quad r + s = \alpha_{33}  & \quad s-r_1 = \alpha_{ 23 }  &   \\
70 ) &  [r_1, r, s, s_1] = [ 2 &  13 &  25 &  31 ] & \quad r + s = \alpha_{33}  &   & \quad r-r_1 = \alpha_{ 9 }  \\
71 ) &  [r_1, r, s, s_1] = [ 2 &  17 &  22 &  31 ] & \quad r + s = \alpha_{33}  & \quad s-r_1 = \alpha_{ 19 }  &   \\
72 ) &  [r_1, r, s, s_1] = [ 2 &  18 &  21 &  31 ] & \quad r + s = \alpha_{33}  &   & \quad r-r_1 = \alpha_{ 14 }  \\
73 ) &  [r_1, r, s, s_1] = [ 3 &  5 &  9 &  10 ] & \quad r + s = \alpha_{14}  & \quad s-r_1 = \alpha_{ 4 }  &   \\
74 ) &  [r_1, r, s, s_1] = [ 3 &  5 &  14 &  15 ] & \quad r + s = \alpha_{19}  & \quad s-r_1 = \alpha_{ 10 }  &   \\
75 ) &  [r_1, r, s, s_1] = [ 3 &  8 &  23 &  26 ] & \quad r + s = \alpha_{29}  &   & \quad r-r_1 = \alpha_{ 2 }  \\
76 ) &  [r_1, r, s, s_1] = [ 3 &  9 &  22 &  26 ] & \quad r + s = \alpha_{29}  &   & \quad r-r_1 = \alpha_{ 4 }  \\
77 ) &  [r_1, r, s, s_1] = [ 3 &  13 &  19 &  26 ] & \quad r + s = \alpha_{29}  & \quad s-r_1 = \alpha_{ 15 }  &   \\
78 ) &  [r_1, r, s, s_1] = [ 3 &  14 &  18 &  26 ] & \quad r + s = \alpha_{29}  &   & \quad r-r_1 = \alpha_{ 10 }  \\
79 ) &  [r_1, r, s, s_1] = [ 4 &  9 &  15 &  19 ] & \quad r + s = \alpha_{23}  &   & \quad r-r_1 = \alpha_{ 3 }  \\
80 ) &  [r_1, r, s, s_1] = [ 4 &  10 &  14 &  19 ] & \quad r + s = \alpha_{23}  &   & \quad r-r_1 = \alpha_{ 5 }  \\
\end{array} $
\caption{ The quartets in $B_6$. }
\label{table_B6_quartets}
\end{table}

\begin{table}[!ht]
$ \qquad
  \tiny
  \begin{array}{llllllll}
1 ) &  [r_1, r, s, s_1] = [ 0 &  2 &  6 &  7 ] & \quad r + s = \alpha_{11}  & \quad s-r_1 = \alpha_{ 1 }  &   \\
2 ) &  [r_1, r, s, s_1] = [ 0 &  3 &  11 &  12 ] & \quad r + s = \alpha_{16}  & \quad s-r_1 = \alpha_{ 7 }  &   \\
3 ) &  [r_1, r, s, s_1] = [ 0 &  6 &  8 &  12 ] & \quad r + s = \alpha_{16}  &   & \quad r-r_1 = \alpha_{ 1 }  \\
4 ) &  [r_1, r, s, s_1] = [ 0 &  4 &  16 &  17 ] & \quad r + s = \alpha_{20}  & \quad s-r_1 = \alpha_{ 12 }  &   \\
5 ) &  [r_1, r, s, s_1] = [ 0 &  6 &  13 &  17 ] & \quad r + s = \alpha_{20}  &   & \quad r-r_1 = \alpha_{ 1 }  \\
6 ) &  [r_1, r, s, s_1] = [ 0 &  9 &  11 &  17 ] & \quad r + s = \alpha_{20}  & \quad s-r_1 = \alpha_{ 7 }  &   \\
7 ) &  [r_1, r, s, s_1] = [ 0 &  5 &  20 &  21 ] & \quad r + s = \alpha_{24}  & \quad s-r_1 = \alpha_{ 17 }  &   \\
8 ) &  [r_1, r, s, s_1] = [ 0 &  6 &  18 &  21 ] & \quad r + s = \alpha_{24}  &   & \quad r-r_1 = \alpha_{ 1 }  \\
9 ) &  [r_1, r, s, s_1] = [ 0 &  10 &  16 &  21 ] & \quad r + s = \alpha_{24}  & \quad s-r_1 = \alpha_{ 12 }  &   \\
10 ) &  [r_1, r, s, s_1] = [ 0 &  11 &  14 &  21 ] & \quad r + s = \alpha_{24}  &   & \quad r-r_1 = \alpha_{ 7 }  \\
11 ) &  [r_1, r, s, s_1] = [ 0 &  4 &  24 &  25 ] & \quad r + s = \alpha_{27}  & \quad s-r_1 = \alpha_{ 21 }  &   \\
12 ) &  [r_1, r, s, s_1] = [ 0 &  6 &  22 &  25 ] & \quad r + s = \alpha_{27}  &   & \quad r-r_1 = \alpha_{ 1 }  \\
13 ) &  [r_1, r, s, s_1] = [ 0 &  10 &  20 &  25 ] & \quad r + s = \alpha_{27}  & \quad s-r_1 = \alpha_{ 17 }  &   \\
14 ) &  [r_1, r, s, s_1] = [ 0 &  11 &  19 &  25 ] & \quad r + s = \alpha_{27}  &   & \quad r-r_1 = \alpha_{ 7 }  \\
15 ) &  [r_1, r, s, s_1] = [ 0 &  15 &  16 &  25 ] & \quad r + s = \alpha_{27}  & \quad s-r_1 = \alpha_{ 12 }  &   \\
16 ) &  [r_1, r, s, s_1] = [ 0 &  3 &  27 &  28 ] & \quad r + s = \alpha_{30}  & \quad s-r_1 = \alpha_{ 25 }  &   \\
17 ) &  [r_1, r, s, s_1] = [ 0 &  6 &  26 &  28 ] & \quad r + s = \alpha_{30}  &   & \quad r-r_1 = \alpha_{ 1 }  \\
18 ) &  [r_1, r, s, s_1] = [ 0 &  9 &  24 &  28 ] & \quad r + s = \alpha_{30}  & \quad s-r_1 = \alpha_{ 21 }  &   \\
19 ) &  [r_1, r, s, s_1] = [ 0 &  11 &  23 &  28 ] & \quad r + s = \alpha_{30}  &   & \quad r-r_1 = \alpha_{ 7 }  \\
20 ) &  [r_1, r, s, s_1] = [ 0 &  14 &  20 &  28 ] & \quad r + s = \alpha_{30}  & \quad s-r_1 = \alpha_{ 17 }  &   \\
21 ) &  [r_1, r, s, s_1] = [ 0 &  16 &  19 &  28 ] & \quad r + s = \alpha_{30}  &   & \quad r-r_1 = \alpha_{ 12 }  \\
22 ) &  [r_1, r, s, s_1] = [ 0 &  2 &  30 &  31 ] & \quad r + s = \alpha_{32}  & \quad s-r_1 = \alpha_{ 28 }  &   \\
23 ) &  [r_1, r, s, s_1] = [ 0 &  6 &  29 &  31 ] & \quad r + s = \alpha_{32}  &   & \quad r-r_1 = \alpha_{ 1 }  \\
24 ) &  [r_1, r, s, s_1] = [ 0 &  8 &  27 &  31 ] & \quad r + s = \alpha_{32}  & \quad s-r_1 = \alpha_{ 25 }  &   \\
25 ) &  [r_1, r, s, s_1] = [ 0 &  11 &  26 &  31 ] & \quad r + s = \alpha_{32}  &   & \quad r-r_1 = \alpha_{ 7 }  \\
26 ) &  [r_1, r, s, s_1] = [ 0 &  13 &  24 &  31 ] & \quad r + s = \alpha_{32}  & \quad s-r_1 = \alpha_{ 21 }  &   \\
27 ) &  [r_1, r, s, s_1] = [ 0 &  16 &  22 &  31 ] & \quad r + s = \alpha_{32}  &   & \quad r-r_1 = \alpha_{ 12 }  \\
28 ) &  [r_1, r, s, s_1] = [ 0 &  18 &  20 &  31 ] & \quad r + s = \alpha_{32}  & \quad s-r_1 = \alpha_{ 17 }  &   \\
29 ) &  [r_1, r, s, s_1] = [ 0 &  1 &  32 &  33 ] & \quad r + s = \alpha_{34}  & \quad s-r_1 = \alpha_{ 31 }  &   \\
30 ) &  [r_1, r, s, s_1] = [ 0 &  6 &  31 &  33 ] & \quad r + s = \alpha_{34}  &   & \quad r-r_1 = \alpha_{ 1 }  \\
31 ) &  [r_1, r, s, s_1] = [ 0 &  7 &  30 &  33 ] & \quad r + s = \alpha_{34}  & \quad s-r_1 = \alpha_{ 28 }  &   \\
32 ) &  [r_1, r, s, s_1] = [ 0 &  11 &  28 &  33 ] & \quad r + s = \alpha_{34}  &   & \quad r-r_1 = \alpha_{ 7 }  \\
33 ) &  [r_1, r, s, s_1] = [ 0 &  12 &  27 &  33 ] & \quad r + s = \alpha_{34}  & \quad s-r_1 = \alpha_{ 25 }  &   \\
34 ) &  [r_1, r, s, s_1] = [ 0 &  16 &  25 &  33 ] & \quad r + s = \alpha_{34}  &   & \quad r-r_1 = \alpha_{ 12 }  \\
35 ) &  [r_1, r, s, s_1] = [ 0 &  17 &  24 &  33 ] & \quad r + s = \alpha_{34}  & \quad s-r_1 = \alpha_{ 21 }  &   \\
36 ) &  [r_1, r, s, s_1] = [ 0 &  20 &  21 &  33 ] & \quad r + s = \alpha_{34}  &   & \quad r-r_1 = \alpha_{ 17 }  \\
37 ) &  [r_1, r, s, s_1] = [ 0 &  6 &  32 &  34 ] & \quad r + s = \alpha_{35}  & \quad s-r_1 = \alpha_{ 31 }  & \quad r-r_1 = \alpha_{ 1 }  \\
38 ) &  [r_1, r, s, s_1] = [ 0 &  11 &  30 &  34 ] & \quad r + s = \alpha_{35}  & \quad s-r_1 = \alpha_{ 28 }  & \quad r-r_1 = \alpha_{ 7 }  \\
39 ) &  [r_1, r, s, s_1] = [ 0 &  16 &  27 &  34 ] & \quad r + s = \alpha_{35}  & \quad s-r_1 = \alpha_{ 25 }  & \quad r-r_1 = \alpha_{ 12 }  \\
40 ) &  [r_1, r, s, s_1] = [ 0 &  20 &  24 &  34 ] & \quad r + s = \alpha_{35}  & \quad s-r_1 = \alpha_{ 21 }  & \quad r-r_1 = \alpha_{ 17 }  \\
41 ) &  [r_1, r, s, s_1] = [ 1 &  3 &  7 &  8 ] & \quad r + s = \alpha_{12}  & \quad s-r_1 = \alpha_{ 2 }  &   \\
42 ) &  [r_1, r, s, s_1] = [ 1 &  4 &  12 &  13 ] & \quad r + s = \alpha_{17}  & \quad s-r_1 = \alpha_{ 8 }  &   \\
43 ) &  [r_1, r, s, s_1] = [ 1 &  7 &  9 &  13 ] & \quad r + s = \alpha_{17}  &   & \quad r-r_1 = \alpha_{ 2 }  \\
44 ) &  [r_1, r, s, s_1] = [ 1 &  5 &  17 &  18 ] & \quad r + s = \alpha_{21}  & \quad s-r_1 = \alpha_{ 13 }  &   \\
45 ) &  [r_1, r, s, s_1] = [ 1 &  7 &  14 &  18 ] & \quad r + s = \alpha_{21}  &   & \quad r-r_1 = \alpha_{ 2 }  \\
46 ) &  [r_1, r, s, s_1] = [ 1 &  10 &  12 &  18 ] & \quad r + s = \alpha_{21}  & \quad s-r_1 = \alpha_{ 8 }  &   \\
47 ) &  [r_1, r, s, s_1] = [ 1 &  4 &  21 &  22 ] & \quad r + s = \alpha_{25}  & \quad s-r_1 = \alpha_{ 18 }  &   \\
48 ) &  [r_1, r, s, s_1] = [ 1 &  7 &  19 &  22 ] & \quad r + s = \alpha_{25}  &   & \quad r-r_1 = \alpha_{ 2 }  \\
49 ) &  [r_1, r, s, s_1] = [ 1 &  10 &  17 &  22 ] & \quad r + s = \alpha_{25}  & \quad s-r_1 = \alpha_{ 13 }  &   \\
50 ) &  [r_1, r, s, s_1] = [ 1 &  12 &  15 &  22 ] & \quad r + s = \alpha_{25}  &   & \quad r-r_1 = \alpha_{ 8 }  \\
51 ) &  [r_1, r, s, s_1] = [ 1 &  3 &  25 &  26 ] & \quad r + s = \alpha_{28}  & \quad s-r_1 = \alpha_{ 22 }  &   \\
52 ) &  [r_1, r, s, s_1] = [ 1 &  7 &  23 &  26 ] & \quad r + s = \alpha_{28}  &   & \quad r-r_1 = \alpha_{ 2 }  \\
53 ) &  [r_1, r, s, s_1] = [ 1 &  9 &  21 &  26 ] & \quad r + s = \alpha_{28}  & \quad s-r_1 = \alpha_{ 18 }  &   \\
54 ) &  [r_1, r, s, s_1] = [ 1 &  12 &  19 &  26 ] & \quad r + s = \alpha_{28}  &   & \quad r-r_1 = \alpha_{ 8 }  \\
55 ) &  [r_1, r, s, s_1] = [ 1 &  14 &  17 &  26 ] & \quad r + s = \alpha_{28}  & \quad s-r_1 = \alpha_{ 13 }  &   \\
56 ) &  [r_1, r, s, s_1] = [ 1 &  2 &  28 &  29 ] & \quad r + s = \alpha_{31}  & \quad s-r_1 = \alpha_{ 26 }  &   \\
57 ) &  [r_1, r, s, s_1] = [ 1 &  7 &  26 &  29 ] & \quad r + s = \alpha_{31}  &   & \quad r-r_1 = \alpha_{ 2 }  \\
58 ) &  [r_1, r, s, s_1] = [ 1 &  8 &  25 &  29 ] & \quad r + s = \alpha_{31}  & \quad s-r_1 = \alpha_{ 22 }  &   \\
59 ) &  [r_1, r, s, s_1] = [ 1 &  12 &  22 &  29 ] & \quad r + s = \alpha_{31}  &   & \quad r-r_1 = \alpha_{ 8 }  \\
60 ) &  [r_1, r, s, s_1] = [ 1 &  13 &  21 &  29 ] & \quad r + s = \alpha_{31}  & \quad s-r_1 = \alpha_{ 18 }  &   \\
61 ) &  [r_1, r, s, s_1] = [ 1 &  17 &  18 &  29 ] & \quad r + s = \alpha_{31}  &   & \quad r-r_1 = \alpha_{ 13 }  \\
62 ) &  [r_1, r, s, s_1] = [ 1 &  7 &  28 &  31 ] & \quad r + s = \alpha_{33}  & \quad s-r_1 = \alpha_{ 26 }  & \quad r-r_1 = \alpha_{ 2 }  \\
63 ) &  [r_1, r, s, s_1] = [ 1 &  12 &  25 &  31 ] & \quad r + s = \alpha_{33}  & \quad s-r_1 = \alpha_{ 22 }  & \quad r-r_1 = \alpha_{ 8 }  \\
64 ) &  [r_1, r, s, s_1] = [ 1 &  17 &  21 &  31 ] & \quad r + s = \alpha_{33}  & \quad s-r_1 = \alpha_{ 18 }  & \quad r-r_1 = \alpha_{ 13 }  \\
65 ) &  [r_1, r, s, s_1] = [ 2 &  4 &  8 &  9 ] & \quad r + s = \alpha_{13}  & \quad s-r_1 = \alpha_{ 3 }  &   \\
66 ) &  [r_1, r, s, s_1] = [ 2 &  5 &  13 &  14 ] & \quad r + s = \alpha_{18}  & \quad s-r_1 = \alpha_{ 9 }  &   \\
67 ) &  [r_1, r, s, s_1] = [ 2 &  8 &  10 &  14 ] & \quad r + s = \alpha_{18}  &   & \quad r-r_1 = \alpha_{ 3 }  \\
68 ) &  [r_1, r, s, s_1] = [ 2 &  4 &  18 &  19 ] & \quad r + s = \alpha_{22}  & \quad s-r_1 = \alpha_{ 14 }  &   \\
69 ) &  [r_1, r, s, s_1] = [ 2 &  8 &  15 &  19 ] & \quad r + s = \alpha_{22}  &   & \quad r-r_1 = \alpha_{ 3 }  \\
70 ) &  [r_1, r, s, s_1] = [ 2 &  10 &  13 &  19 ] & \quad r + s = \alpha_{22}  & \quad s-r_1 = \alpha_{ 9 }  &   \\
71 ) &  [r_1, r, s, s_1] = [ 2 &  3 &  22 &  23 ] & \quad r + s = \alpha_{26}  & \quad s-r_1 = \alpha_{ 19 }  &   \\
72 ) &  [r_1, r, s, s_1] = [ 2 &  8 &  19 &  23 ] & \quad r + s = \alpha_{26}  &   & \quad r-r_1 = \alpha_{ 3 }  \\
73 ) &  [r_1, r, s, s_1] = [ 2 &  9 &  18 &  23 ] & \quad r + s = \alpha_{26}  & \quad s-r_1 = \alpha_{ 14 }  &   \\
74 ) &  [r_1, r, s, s_1] = [ 2 &  13 &  14 &  23 ] & \quad r + s = \alpha_{26}  &   & \quad r-r_1 = \alpha_{ 9 }  \\
75 ) &  [r_1, r, s, s_1] = [ 2 &  8 &  22 &  26 ] & \quad r + s = \alpha_{29}  & \quad s-r_1 = \alpha_{ 19 }  & \quad r-r_1 = \alpha_{ 3 }  \\
76 ) &  [r_1, r, s, s_1] = [ 2 &  13 &  18 &  26 ] & \quad r + s = \alpha_{29}  & \quad s-r_1 = \alpha_{ 14 }  & \quad r-r_1 = \alpha_{ 9 }  \\
77 ) &  [r_1, r, s, s_1] = [ 3 &  5 &  9 &  10 ] & \quad r + s = \alpha_{14}  & \quad s-r_1 = \alpha_{ 4 }  &   \\
78 ) &  [r_1, r, s, s_1] = [ 3 &  4 &  14 &  15 ] & \quad r + s = \alpha_{19}  & \quad s-r_1 = \alpha_{ 10 }  &   \\
79 ) &  [r_1, r, s, s_1] = [ 3 &  9 &  10 &  15 ] & \quad r + s = \alpha_{19}  &   & \quad r-r_1 = \alpha_{ 4 }  \\
80 ) &  [r_1, r, s, s_1] = [ 3 &  9 &  14 &  19 ] & \quad r + s = \alpha_{23}  & \quad s-r_1 = \alpha_{ 10 }  & \quad r-r_1 = \alpha_{ 4 }  \\
  \end{array} $
\caption{ The quartets in $C_6$. }
\label{table_C6_quartets}
\end{table}

\begin{table}[!ht]
$ \qquad
  \begin{array}{llllllll}
1 ) &  [r_1, r, s, s_1] = [ 0 &  2 &  4 &  5 ] & \quad r + s = \alpha_{7}  & \quad s-r_1 = \alpha_{ 1 }  &   \\
2 ) &  [r_1, r, s, s_1] = [ 0 &  2 &  7 &  8 ] & \quad r + s = \alpha_{10}  & \quad s-r_1 = \alpha_{ 5 }  &   \\
3 ) &  [r_1, r, s, s_1] = [ 0 &  3 &  7 &  9 ] & \quad r + s = \alpha_{11}  & \quad s-r_1 = \alpha_{ 5 }  &   \\
4 ) &  [r_1, r, s, s_1] = [ 0 &  4 &  6 &  9 ] & \quad r + s = \alpha_{11}  &   & \quad r-r_1 = \alpha_{ 1 }  \\
5 ) &  [r_1, r, s, s_1] = [ 0 &  2 &  11 &  12 ] & \quad r + s = \alpha_{14}  & \quad s-r_1 = \alpha_{ 9 }  &   \\
6 ) &  [r_1, r, s, s_1] = [ 0 &  3 &  10 &  12 ] & \quad r + s = \alpha_{14}  & \quad s-r_1 = \alpha_{ 8 }  &   \\
7 ) &  [r_1, r, s, s_1] = [ 0 &  6 &  7 &  12 ] & \quad r + s = \alpha_{14}  & \quad s-r_1 = \alpha_{ 5 }  &   \\
8 ) &  [r_1, r, s, s_1] = [ 0 &  3 &  14 &  15 ] & \quad r + s = \alpha_{17}  & \quad s-r_1 = \alpha_{ 12 }  &   \\
9 ) &  [r_1, r, s, s_1] = [ 0 &  6 &  11 &  15 ] & \quad r + s = \alpha_{17}  & \quad s-r_1 = \alpha_{ 9 }  &   \\
10 ) &  [r_1, r, s, s_1] = [ 0 &  4 &  21 &  22 ] & \quad r + s = \alpha_{23}  &   & \quad r-r_1 = \alpha_{ 1 }  \\
11 ) &  [r_1, r, s, s_1] = [ 0 &  7 &  20 &  22 ] & \quad r + s = \alpha_{23}  &   & \quad r-r_1 = \alpha_{ 5 }  \\
12 ) &  [r_1, r, s, s_1] = [ 0 &  10 &  19 &  22 ] & \quad r + s = \alpha_{23}  &   & \quad r-r_1 = \alpha_{ 8 }  \\
13 ) &  [r_1, r, s, s_1] = [ 0 &  11 &  18 &  22 ] & \quad r + s = \alpha_{23}  &   & \quad r-r_1 = \alpha_{ 9 }  \\
14 ) &  [r_1, r, s, s_1] = [ 0 &  13 &  17 &  22 ] & \quad r + s = \alpha_{23}  & \quad s-r_1 = \alpha_{ 15 }  &   \\
15 ) &  [r_1, r, s, s_1] = [ 0 &  14 &  16 &  22 ] & \quad r + s = \alpha_{23}  &   & \quad r-r_1 = \alpha_{ 12 }  \\
16 ) &  [r_1, r, s, s_1] = [ 1 &  3 &  5 &  6 ] & \quad r + s = \alpha_{9}  & \quad s-r_1 = \alpha_{ 2 }  &   \\
17 ) &  [r_1, r, s, s_1] = [ 1 &  4 &  8 &  10 ] & \quad r + s = \alpha_{13}  &   & \quad r-r_1 = \alpha_{ 0 }  \\
18 ) &  [r_1, r, s, s_1] = [ 1 &  5 &  7 &  10 ] & \quad r + s = \alpha_{13}  &   & \quad r-r_1 = \alpha_{ 2 }  \\
19 ) &  [r_1, r, s, s_1] = [ 1 &  3 &  13 &  14 ] & \quad r + s = \alpha_{16}  & \quad s-r_1 = \alpha_{ 10 }  &   \\
20 ) &  [r_1, r, s, s_1] = [ 1 &  4 &  12 &  14 ] & \quad r + s = \alpha_{16}  &   & \quad r-r_1 = \alpha_{ 0 }  \\
21 ) &  [r_1, r, s, s_1] = [ 1 &  5 &  11 &  14 ] & \quad r + s = \alpha_{16}  &   & \quad r-r_1 = \alpha_{ 2 }  \\
22 ) &  [r_1, r, s, s_1] = [ 1 &  7 &  9 &  14 ] & \quad r + s = \alpha_{16}  & \quad s-r_1 = \alpha_{ 6 }  &   \\
23 ) &  [r_1, r, s, s_1] = [ 1 &  3 &  16 &  17 ] & \quad r + s = \alpha_{19}  & \quad s-r_1 = \alpha_{ 14 }  &   \\
24 ) &  [r_1, r, s, s_1] = [ 1 &  4 &  15 &  17 ] & \quad r + s = \alpha_{19}  &   & \quad r-r_1 = \alpha_{ 0 }  \\
25 ) &  [r_1, r, s, s_1] = [ 1 &  9 &  11 &  17 ] & \quad r + s = \alpha_{19}  &   & \quad r-r_1 = \alpha_{ 6 }  \\
26 ) &  [r_1, r, s, s_1] = [ 1 &  5 &  20 &  21 ] & \quad r + s = \alpha_{22}  &   & \quad r-r_1 = \alpha_{ 2 }  \\
27 ) &  [r_1, r, s, s_1] = [ 1 &  8 &  19 &  21 ] & \quad r + s = \alpha_{22}  & \quad s-r_1 = \alpha_{ 17 }  &   \\
28 ) &  [r_1, r, s, s_1] = [ 1 &  9 &  18 &  21 ] & \quad r + s = \alpha_{22}  &   & \quad r-r_1 = \alpha_{ 6 }  \\
29 ) &  [r_1, r, s, s_1] = [ 1 &  12 &  16 &  21 ] & \quad r + s = \alpha_{22}  & \quad s-r_1 = \alpha_{ 14 }  &   \\
30 ) &  [r_1, r, s, s_1] = [ 1 &  13 &  15 &  21 ] & \quad r + s = \alpha_{22}  &   & \quad r-r_1 = \alpha_{ 10 }  \\
31 ) &  [r_1, r, s, s_1] = [ 2 &  3 &  8 &  9 ] & \quad r + s = \alpha_{12}  & \quad s-r_1 = \alpha_{ 5 }  &   \\
32 ) &  [r_1, r, s, s_1] = [ 2 &  5 &  6 &  9 ] & \quad r + s = \alpha_{12}  & \quad s-r_1 = \alpha_{ 3 }  & \quad r-r_1 = \alpha_{ 1 }  \\
33 ) &  [r_1, r, s, s_1] = [ 2 &  5 &  14 &  16 ] & \quad r + s = \alpha_{18}  & \quad s-r_1 = \alpha_{ 11 }  & \quad r-r_1 = \alpha_{ 1 }  \\
34 ) &  [r_1, r, s, s_1] = [ 2 &  6 &  13 &  16 ] & \quad r + s = \alpha_{18}  &   & \quad r-r_1 = \alpha_{ 3 }  \\
35 ) &  [r_1, r, s, s_1] = [ 2 &  7 &  12 &  16 ] & \quad r + s = \alpha_{18}  & \quad s-r_1 = \alpha_{ 9 }  & \quad r-r_1 = \alpha_{ 4 }  \\
36 ) &  [r_1, r, s, s_1] = [ 2 &  8 &  11 &  16 ] & \quad r + s = \alpha_{18}  &   & \quad r-r_1 = \alpha_{ 5 }  \\
37 ) &  [r_1, r, s, s_1] = [ 2 &  9 &  10 &  16 ] & \quad r + s = \alpha_{18}  & \quad s-r_1 = \alpha_{ 7 }  &   \\
38 ) &  [r_1, r, s, s_1] = [ 2 &  3 &  18 &  19 ] & \quad r + s = \alpha_{20}  & \quad s-r_1 = \alpha_{ 16 }  &   \\
39 ) &  [r_1, r, s, s_1] = [ 2 &  5 &  17 &  19 ] & \quad r + s = \alpha_{20}  &   & \quad r-r_1 = \alpha_{ 1 }  \\
40 ) &  [r_1, r, s, s_1] = [ 2 &  6 &  16 &  19 ] & \quad r + s = \alpha_{20}  &   & \quad r-r_1 = \alpha_{ 3 }  \\
41 ) &  [r_1, r, s, s_1] = [ 2 &  7 &  15 &  19 ] & \quad r + s = \alpha_{20}  &   & \quad r-r_1 = \alpha_{ 4 }  \\
42 ) &  [r_1, r, s, s_1] = [ 2 &  9 &  14 &  19 ] & \quad r + s = \alpha_{20}  & \quad s-r_1 = \alpha_{ 11 }  &   \\
43 ) &  [r_1, r, s, s_1] = [ 2 &  11 &  12 &  19 ] & \quad r + s = \alpha_{20}  & \quad s-r_1 = \alpha_{ 9 }  &   \\
44 ) &  [r_1, r, s, s_1] = [ 2 &  6 &  18 &  20 ] & \quad r + s = \alpha_{21}  & \quad s-r_1 = \alpha_{ 16 }  & \quad r-r_1 = \alpha_{ 3 }  \\
45 ) &  [r_1, r, s, s_1] = [ 2 &  8 &  17 &  20 ] & \quad r + s = \alpha_{21}  &   & \quad r-r_1 = \alpha_{ 5 }  \\
46 ) &  [r_1, r, s, s_1] = [ 2 &  10 &  15 &  20 ] & \quad r + s = \alpha_{21}  &   & \quad r-r_1 = \alpha_{ 7 }  \\
47 ) &  [r_1, r, s, s_1] = [ 2 &  12 &  14 &  20 ] & \quad r + s = \alpha_{21}  & \quad s-r_1 = \alpha_{ 11 }  & \quad r-r_1 = \alpha_{ 9 }  \\
48 ) &  [r_1, r, s, s_1] = [ 3 &  6 &  9 &  12 ] & \quad r + s = \alpha_{15}  & \quad s-r_1 = \alpha_{ 5 }  & \quad r-r_1 = \alpha_{ 2 }  \\
  \end{array} $
\medskip
\caption{ The quartets in $F_4$. The first $30$ lines are the simple quartets with
$\abs{r_1} = \sqrt{2}$. }
\label{table_F4_quartets}
\end{table}

\end{appendix}

\end{document}